\documentclass[reqno]{amsart} 
 
\usepackage{hyperref}
\usepackage{amsmath}
\usepackage{amssymb}
\usepackage{amsfonts}
\usepackage{amsthm}
\usepackage{latexsym}
\usepackage{esint}
\usepackage{mathtools}
\usepackage{mathrsfs}
\usepackage{graphicx}
\usepackage{wrapfig}
\usepackage{bbm}
\usepackage{color}
\usepackage{bm}
\usepackage[top=1in, bottom=1.25in, left=1.10in, right=1.10in]{geometry}
\usepackage{todonotes}
\allowdisplaybreaks

\newtheorem{theorem}{Theorem}[section]

\newtheorem{proposition}[theorem]{Proposition}
\newtheorem{corollary}[theorem]{Corollary}

\theoremstyle{definition}
\newtheorem{definition}[theorem]{Definition}

\theoremstyle{remark}
\newtheorem{remark}[theorem]{Remark}

\numberwithin{equation}{section}

\newcommand{\bx}{\mathbf{x}}
\newcommand{\by}{{y}}

\newcommand{\bq}{\mathbf{q}}
\newcommand{\bu}{\mathbf{u}}
\newcommand{\bw}{\mathbf{w}}
\newcommand{\bfu}{\mathbf{u}}

\newcommand{\bff}{\mathbf{f}}

\newcommand{\bv}{\mathbf{v}}
\newcommand{\bn}{\mathbf{n}}

\newcommand{\dy}{\, \mathrm{d}y}
\newcommand{\dd}{\,\mathrm{d}}
\newcommand{\dq}{\, \mathrm{d} \mathbf{q}}
\newcommand{\dx}{\, \mathrm{d} \mathbf{x}}

\newcommand{\dt}{\, \mathrm{d}t}
\newcommand{\ds}{\, \mathrm{d}s}

\newcommand{\divx}{\mathrm{div}_{\mathbf{x}}}
\newcommand{\divq}{\mathrm{div}_{\mathbf{q}}}

\newcommand{\delx}{\Delta_{\mathbf{x}}}

\newcommand{\nabx}{\nabla_{\mathbf{x}}}
\newcommand{\naby}{\partial_y}
\newcommand{\nabq}{\nabla_{\mathbf{q}}}

\newcommand{\Delx}{\Delta_{\mathbf{x}}}

\newcommand{\bT}{\mathbb{T}}
\newcommand{\bU}{\mathbb{U}}

\newcommand{\R}{\mathbb{R}}

\newcommand{\Oeta}{\Omega_{\eta}}
\newcommand{\Ozeta}{\Omega_{\zeta}}

\begin{document} 
%\title[The elastic, corotational, and vanishing center-of-mass limit]{The elastic, corotational, and vanishing center-of-mass limit of polymeric fluid-structure interaction of {O}ldroyd-{B} type}

%\title[The elastic, corotational, and vanishing center-of-mass limit]{The tripartite limit of polymeric fluid-structure interaction of {O}ldroyd-{B} type} 

\title[The vanishing center-of-mass limit]{Vanishing center-of-mass limit of the 2{D}-1{D} corotational {O}ldroyd-{B} polymeric fluid-structure interaction problem}

%%    Information for first author
%\author{Dominic Breit}
%%    Address of record for the research reported here
%\address{Department of Mathematics, Heriot-Watt University, Edinburgh, EH14 4AS, United Kingdom}
%\email{d.breit@hw.ac.uk}
%
%\thanks{The authors would like to thank S. Schwarzacher and E. S\"uli for valuable suggestions.}
\author{Prince Romeo Mensah}
\address{Institut f\"ur Mathematik,
Technische Universit\"at Clausthal,
Erzstra{\ss}e 1,
38678 Clausthal-Zellerfeld}
\email{prince.romeo.mensah@tu-clausthal.de \\ orcid ID: 0000-0003-4086-2708}
%%\thanks{Support information for the second author.}

%    General info
\subjclass[2020]{76D03; 74F10; 35Q30; 35Q84; 82D60}

\date{\today}


\keywords{Incompressible Navier--Stokes--Fokker--Planck system, Oldroyd-B, Fluid-Structure interaction, Polymeric fluid, Corotational fluid}

\begin{abstract} 
We consider the Oldroyd-B model for a  two-dimensional dilute corotational polymer fluid with center-of-mass diffusion that is interacting with a one-dimensional viscoelastic shell. We show that any family of strong solutions of the system described above that is parametrized by the center-of-mass diffusion coefficient converges, as the coefficient goes to zero, to a weak solution of a  corotational polymer fluid-structure interaction system without center-of-mass diffusion but with essentially bounded polymer number density and extra stress tensor. 
As a consequence, we also obtain a weak-strong uniqueness result that says that the weak solution of the latter is unique in the class of the strong solution of the former as the center-of-mass diffusion vanishes.

%
%*************
% The elastic, corotational, and vanishing center-of-mass limit of polymeric fluid-structure interaction of {O}ldroyd-{B} type
\end{abstract}

\maketitle
%\tableofcontents

\section{Introduction}
The analysis of the interaction of polymeric fluids with flexible structures was recently initiated in \cite{breit2021incompressible}  leading to a connection between two thriving fields in continuum mechanics, i.e. polymeric fluid analysis and fluid-structure interaction problems.  A key challenge to merging these two fields stemmed from their modeling. The typical analysis of a classical fluid-structure interaction problem involves the interplay between the macroscopic description of a solvent and a structure whereas polymeric fluids analysis deals with the interaction of the macroscopic description of a solvent with a solute on either the microscopic, mesoscopic, or macroscopic scale.  Furthermore, the high molecular weight of a polymeric fluid (as compared to a Newtonian or Ideal fluid, say) makes it comparatively solid-like and thus affects the surface force at the boundary. These inhibiting issues have since been resolved in \cite{breit2021incompressible} opening the field to further research. 
 
In this work, we are interested in the interaction of \textit{corotational} polymeric fluids with flexible structures. A corotational polymer fluid is a polymer fluid whose drag term for the equation of the solute is driven by the vorticity tensor of the solvent. The interaction of vorticity with the velocity gradients present in the modeling of the solvent component of the polymer solution is a key ingredient in the creation and maintenance of turbulence. The rigorous understanding of turbulence in simple Newtonian fluids is already a challenging task and it is not the focus of this work. Instead, we focus on these chaotic corotational polymer fluids within viscoelastic structures and the interaction between these polymer fluids and the structures. To properly motivate this work, we begin with the system setup.

\subsection{Setup}
We consider a spatial reference domain  $\Omega \subset \mathbb{R}^2$ whose boundary $\partial\Omega$ may consist of a flexible part $\omega \subset \mathbb{R}$ and a rigid part $\Gamma \subset\mathbb{R}$. However, because the analysis at the rigid part is significantly simpler, we shall identify the whole of $\partial \Omega$ with $\omega$ and endow it with periodic boundary conditions. Let $I:=(0,T)$ represent a time interval for a given constant $T>0$. 
We represent the time-dependent  displacement of the structure by $\eta:\overline{I}\times\omega\rightarrow(-L,L)$ where $L>0$ is a fixed length of the tubular neighbourhood of $\partial\Omega$ given by
\begin{align*}
S_L:=\{\bx\in \mathbb{R}^2\,:\, \mathrm{dist}(\bx,\partial\Omega
)<L \}.
\end{align*}
Now, for some $k\in\mathbb{N}$ large enough, we assume that $\partial\Omega$  is parametrized by an injective mapping $\bm{\varphi}\in C^k(\omega;\mathbb{R}^2)$ with $\naby \bm{\varphi}\neq0$ such that
\begin{align*}
\partial{\Omega_{\eta(t)}}=\big\{\bm{\varphi}_{\eta(t)}:=\bm{\varphi}(\by)+\bn(\by)\eta(t,\by)\, :\, t\in I, \by\in \omega\big\}.
\end{align*}
The set $\partial{\Omega_{\eta(t)}}$ represents the boundary of the flexible domain at any instant of time $t\in I$ and the vector $\bn(y)$ is a unit normal at the point $y\in \omega$. 
We also let $\bn_{\eta(t)}(y)$ be the corresponding normal of $\partial{\Omega_{\eta(t)}}$ at the spacetime point $y\in \omega$ and $t\in I$. Then for $L>0$ sufficiently small, $\bn_{\eta(t)}(y)$ is close to $\bn(y)$ and $\bm{\varphi}_{\eta(t)}$ is close to $\bm{\varphi}$. Since $\naby \bm{\varphi}\neq0$,  it will follow that
\begin{align*}
\naby \bm{\varphi}_{\eta(t)} \neq0 \quad\text{ and }\quad \bn(y)\cdot \bn_{\eta(t)}(y)\neq 0 
\end{align*}
for $y\in \omega$ and $t\in I$. Thus, in particular, there is no loss of strict positivity of the Jacobian determinant provided that $\Vert \eta\Vert_{L^\infty(I\times\omega)}<L$.

For the interior points, we  transform the  reference domain $\Omega$ into a time-dependent moving domain $\Omega_{\eta(t)}$  whose state at time $t\in\overline{I}$ is given by
\begin{align*}
\Omega_{\eta(t)}
 =\big\{
 \bm{\Psi}_{\eta(t)}(\bx):\, \bx \in \Omega 
  \big\}.
\end{align*}
Here,
\begin{align*}
\bm{\Psi}_{\eta(t)}(\bx)=
\begin{cases}
\bx+\bn(\by(\bx))\eta(t,\by(\bx))\phi(s(\bx))     & \quad \text{if } \mathrm{dist}(\bx,\partial\Omega)<L,\\
    \bx & \quad \text{elsewhere } 
  \end{cases}
\end{align*}
is the Hanzawa transform with inverse $\bm{\Psi}_{-\eta(t)}$ and where for a point $\bx$ in the neighbourhood of $\partial\Omega$, the vector $\bn(y(\bx))$ is the unit normal at the point $y(\bx)=\mathrm{arg min}_{y\in\omega}\vert\bx -\bm{\varphi}(y)\vert$. Also, $s(\bx)=(\bx-\bm{\varphi}(y(\bx)))\cdot\bn(y(\bx))$ and $\phi\in C^\infty(\mathbb{R})$ is a cut-off function that is $\phi\equiv0$ in the neighbourhood of $-L$ and $\phi\equiv1$ in the neighbourhood of $0$. Note that $\bm{\Psi}_{\eta(t)}(\bx)$ can be rewritten as
\begin{align*}
\bm{\Psi}_{\eta(t)}(\bx)=
\begin{cases}
\bm{\varphi}(y(\bx))+\bn(\by(\bx))[s(\bx)+\eta(t,\by(\bx))\phi(s(\bx)) ]    & \quad \text{if } \mathrm{dist}(\bx,\partial\Omega)<L,\\
    \bx & \quad \text{elsewhere. } 
  \end{cases}
\end{align*}
\\
With the above preparation in hand, we consider the corotational Oldroyd-B model for the flow of a  dilute polymeric fluid interacting with a flexible structure in  the closure of the deformed spacetime cylinder
\begin{align*}
I\times\Oeta:=\bigcup_{t\in I}\{t\}\times\Oeta
\end{align*}
with $\Oeta:={\Omega_{\eta(t)}}$. The unknowns consists of a structure displacement function $\eta:(t, \by)\in I \times \omega \mapsto   \eta(t,\by)\in \mathbb{R}$, a fluid velocity field $\mathbf{u}:(t, \mathbf{x})\in I \times \Oeta \mapsto  \mathbf{u}(t, \mathbf{x}) \in \mathbb{R}^2$, a pressure function $p:(t, \mathbf{x})\in I \times \Oeta \mapsto  p(t, \mathbf{x}) \in \mathbb{R}$, a polymer number density $\rho :(t, \mathbf{x} )\in I \times \Oeta  \mapsto \rho(t, \mathbf{x} ) \in \mathbb{R}$
and an extra stress tensor $\bT :(t, \mathbf{x} )\in I \times \Oeta  \mapsto \bT (t, \mathbf{x} ) \in \mathbb{R}^{2\times2}$
 such that the system of equations 
\begin{align}
\divx \bu=0, \label{divfree}
\\
\partial_t \rho+ (\mathbf{u}\cdot \nabx) \rho
= 
\varepsilon\Delx \rho 
,\label{rhoEqu}
\\
\partial_t \bu  + (\mathbf{u}\cdot \nabx)\mathbf{u} 
= 
 \delx \bu -\nabx p 
+
\divx   \bT, \label{momEq}
\\
 \partial_t^2 \eta - \gamma\partial_t\partial_y^2 \eta +  \partial_y^4 \eta =  - ( \mathbb{S}\bn_\eta )\circ \bm{\varphi}_\eta\cdot\bn \,\det(\partial_y\bm{\varphi}_\eta), 
\label{shellEQ}
\\
\partial_t \bT + (\mathbf{u}\cdot \nabx) \bT
=
\mathbb{W}(\nabx \bu)\bT + \bT\mathbb{W}((\nabx \bu)^\top) - 2(\bT - \rho \mathbb{I})+\varepsilon\Delx \bT \label{solute}
\end{align}
holds on $I\times\Oeta\subset \mathbb R^{1+2}$ (with \eqref{shellEQ} posed on $I\times\omega\subset \mathbb R^{1+1}$) where
\begin{align*}
\mathbb{S}= (\nabx \bu +(\nabx \bu)^\top) -p\mathbb{I}+  \bT,
\qquad
\mathbb{W}(\nabx \bu)=\tfrac{1}{2}(\nabx \bu-(\nabx \bu)^\top).
\end{align*}
The parameters $\varepsilon$ and $\gamma$ are positive constants, $\bn_\eta$ is the normal at $\partial\Oeta$ and $\mathbb{I}$ is the identity matrix.
We complement \eqref{divfree}--\eqref{solute} with the following initial and boundary conditions
\begin{align}
&\eta(0,\cdot)=\eta_0(\cdot), \qquad\partial_t\eta(0,\cdot)=\eta_\star(\cdot) & \text{in }\omega,
\\
&\bu(0,\cdot)=\bu_0(\cdot) & \text{in }\Omega_{\eta_0},
\\
&\rho(0,\cdot)=\rho_0(\cdot),\quad\bT(0,\cdot)=\bT_0(\cdot) &\text{in }\Omega_{\eta_0},
\label{initialCondSolv}
\\
& 
\bn_{\eta}\cdot\nabx\rho=0,\qquad
\bn_{\eta}\cdot\nabx\bT=0 &\text{on }I\times\partial\Omega_{\eta}.
\label{bddCondSolv}
\end{align}
Furthermore, we impose periodicity on the boundary of $\omega$ (with mean-zero elements in $\omega$) and the following interface condition
\begin{align} 
\label{interface}
&\bu\circ\bm{\varphi}_\eta=(\partial_t\eta)\bn & \text{on }I\times \omega
\end{align}
at the flexible part of the boundary with normal $\bn$.
\\
The two unknowns $\rho$ and $\bT$ for the solute component of the polymer fluid are related via the identities
\begin{align*} 
\bT(t, \bx)= \int_{B} f(t,\bx,\bq)\bq\otimes\bq \dq,
\qquad
\rho(t, \bx)= \int_{B} f(t,\bx,\bq) \dq,
\end{align*}
respectively, where for  $B=\mathbb{R}^2$ with elements $\bq\in B$, the function $f$ is the probability density function ($f\geq0$ a.e. on $\overline{I}\times\Oeta\times B$) satisfying the
Fokker--Planck equation
\begin{align}
\label{fokkerPlanck}
\partial_t f+\divx(\bu f)+\divq(\mathbb{W}(\nabx\bu)\bq f)=\varepsilon\Delx f+ \divq(M\nabq(f/M))
\end{align}
in $I\times\Oeta\times B$ for a Hookean dumbbell spring potential and Maxwellian
\begin{align*} 
U\Big(\frac{1}{2}|\bq|^2 \Big)=\frac{1}{2}|\bq|^2, \qquad\qquad M=\frac{\exp(-U(\tfrac{1}{2}|\bq|^2))}{\int_B\exp(-U(\tfrac{1}{2}|\bq|^2))\dq},
\end{align*}
respectively.

As stated earlier, polymeric fluid-structure interaction problems were recently initiated in \cite{breit2021incompressible} where the authors showed the existence of weak solutions to a dilute solute-solvent-structure mutually coupled system. This system consisted of the $3$-D noncorotational Fokker--Planck equation \eqref{fokkerPlanck} (where $\nabx\bu$ replaces $\mathbb{W}(\nabx\bu)$) for the mesoscopic description of the solute, the  $3$-D incompressible Navier--Stokes equation \eqref{divfree} and \eqref{momEq} giving the macroscopic description of the solvent, and with a  $2$-D structure modeled by a shell equation \eqref{shellEQ} of Koiter type (where the term $- \gamma\partial_t\partial_y^2 \eta +  \partial_y^4 \eta$ is replaced by the gradient of the so-called \textit{Koiter energy}). Uniqueness is unknown for this system but the solutions exist until potential degeneracies occur with the Koiter energy or with the structure deformation.
When the  $2$-D Koiter shell in \cite{breit2021incompressible} is replaced by the $2$-D viscoelastic shell equation \eqref{shellEQ}, the extension to the existence of a unique local-in-time strong solution was then shown in \cite{breit2023existence}. 
Note that for fixed spatial domains subject to periodic boundary conditions, one can construct solutions that are spatially more regular \cite{breit2021local}  than strong solutions. The corresponding result for the system with a structure displacement remains an interesting open problem even in lower dimensions.

In this paper, we shall discuss both strong and weak solutions for the dilute $2$-D/$1$-D corotational solute-solvent-structure system \eqref{divfree}--\eqref{interface} of Oldroyd-B type.

\subsection{Concepts of solution and main results}
We begin this subsection with a precise definition of a 
 strong solutions of \eqref{divfree}--\eqref{interface}.
\begin{definition}[Strong solution]
\label{def:strongSolution}
Let $(\rho_0, \bT_0, \bu_0, \eta_0, \eta_\star)$
be a dataset that satisfies
\begin{equation}
\begin{aligned}
\label{mainDataForAllStrong}
%&\bff \in L^2(I;L^{2}_\mathrm{loc}(\mathbb{R}^2)),
%\qquad g\in L^2(I;L^{2}(\omega)) ,
%\\
&
\eta_0 \in W^{3,2}(\omega) \text{ with } \Vert \eta_0 \Vert_{L^\infty( \omega)} < L, \quad \eta_\star \in W^{1,2}(\omega),
\\&\bu_0 \in W^{1,2}_{\divx}(\Omega_{\eta_0} )\text{ is such that }\bu_0 \circ \bm{\varphi}_{\eta_0} =\eta_\star \bn \text{ on } \omega,
\\&
\rho_0\in W^{1,2}(\Omega_{\eta_0}), \quad
\bT_0\in W^{1,2}(\Omega_{\eta_0}),
\\&
\rho_0\geq 0,\,\, \bT_0>0 \quad \text{a.e. in } \Omega_{\eta_0}.
\end{aligned}
\end{equation}
We call 
$(\eta, \bu, p, \rho,\bT)$ 
a \textit{strong solution} of   \eqref{divfree}--\eqref{interface} with dataset $(\rho_0, \bT_0, \bu_0, \eta_0, \eta_\star)$  if:
\begin{itemize}
%\item[(a)]  $(\eta, \bu, \rho,\bT)$ is a weak solution of   \eqref{divfree}--\eqref{interface};
\item the structure-function $\eta$ is such that $
\Vert \eta \Vert_{L^\infty(I \times \omega)} <L$ and
\begin{align*}
\eta \in &W^{1,\infty}\big(I;W^{1,2}(\omega)  \big)\cap L^{\infty}\big(I;W^{3,2}(\omega)  \big) \cap  W^{1,2}\big(I; W^{2,2}(\omega)  \big)
\\&\cap  W^{2,2}\big(I;L^{2}(\omega)  \big) \cap  L^{2}\big(I;W^{4,2}(\omega)  \big);
\end{align*}
\item the velocity $\bu$ is such that $\bu  \circ \bm{\varphi}_{\eta} =(\partial_t\eta)\bn$ on $I\times \omega$ and
\begin{align*} 
\bu\in  W^{1,2} \big(I; L^2_{\divx}(\Oeta ) \big)\cap L^2\big(I;W^{2,2}(\Oeta)  \big);
\end{align*}
\item the pressure $p$ is such that 
\begin{align*}
p\in L^2\big(I;W^{1,2}(\Oeta)  \big);
\end{align*}
\item the pair $(\rho,\bT)$ is such that 
\begin{align*}
\rho,\bT \in W^{1,2}\big(I;L^{2}(\Oeta)  \big) \cap L^\infty\big(I;W^{1,2}(\Oeta)  \big)\cap L^2\big(I;W^{2,2}(\Oeta)  \big);
\end{align*}
\item equations \eqref{divfree}--\eqref{solute} are satisfied a.e. in spacetime with $\eta(0)=\eta_0$ and $\partial_t\eta(0)=\eta_\star$ a.e. in $\omega$, as well as $\bfu(0)=\bfu_0$, $\rho(0)=\rho_0$ and $\bT(0)=\bT_0$ a.e. in $\Omega_{\eta_0}$.
\end{itemize}
\end{definition} 
The existence of a unique strong solution of  \eqref{divfree}--\eqref{interface} in the sense of Definition \ref{def:strongSolution} has recently been shown in \cite{mensah2023weak} (see also \cite{mensah2024conditionally} for a conditional regularity result). 
In the mathematical analysis of polymeric fluids,  there are strong opinions on whether or not the center-of-mass diffusion terms (the terms with the $\varepsilon$ parameters) in \eqref{divfree}--\eqref{solute} should be added. A school of thought \cite{barrett2017existenceOldroyd, degond2009kinetic, schieber2006generalized} gives justifications for its inclusion while some others \cite{guillope1990existence, lin2005on} ignore it. The rigorous derivation (see for example \cite{suli2018mckean}) from the microscopic scale leads to the center-of-mass diffusion terms in the mesoscopic description \eqref{fokkerPlanck} and thus persist on the macroscopic level \eqref{solute}. However, the importance of this stress diffusion is also known to decrease as the length scale of the problem increases \cite{bhave1991kinetic}. From earlier works \cite{breit2021incompressible, breit2023existence} on the Navier--Stokes--Fokker--Planck equation, this suggests that a strong solution (which typically exists locally in time, cf. \cite{breit2023existence}) may converge to a weak solution (which will ideally exist globally in time barring any degeneracies, cf. \cite{breit2021incompressible}) as the center-of-mass diffusion parameter goes to zero. We however note that the notion of a strong solution, as constructed in \cite{breit2023existence}, only applies in spacetime and that the solution is interpreted weakly in the third independent variable, i.e. concerning the conformation vector. On the other hand, the strong solution of  \eqref{divfree}--\eqref{interface} constructed \cite{mensah2023weak} does not have this mixed character since it is fully macroscopic with just a spacetime dependency. With the above discussion in hand,
we show that the system \eqref{divfree}--\eqref{solute} which models a polymeric fluid of Oldroyd-B type with center-of-mass diffusion that is interacting with a viscoelastic shell converges to the corotational polymeric fluid model of Oldroyd-B type \textit{without} center-of-mass diffusion that is also interacting with a viscoelastic shell. The latter model is described by the following system of equations
\begin{align}
\divx \bv=0, \label{divfreeL}
\\
\partial_t \sigma+ (\bv\cdot \nabx) \sigma
= 
0
,\label{rhoEquL}
\\
\partial_t \bv  + (\bv\cdot \nabx)\bv
= 
  \delx \bv -\nabx \pi 
+ 
\divx   \bU, \label{momEqL}
\\
 \partial_t^2 \zeta - \gamma\partial_t\partial_y^2 \zeta +\partial_y^4 \zeta =  - ( \mathbb{S}_*\bn_\zeta )\circ \bm{\varphi}_\zeta\cdot\bn \,\det(\partial_y\bm{\varphi}_\zeta), 
\label{shellEQL}
\\
\partial_t \bU + (\bv\cdot \nabx) \bU
=
\mathbb{W}(\nabx \bv)\bU + \bU\mathbb{W}((\nabx \bv)^\top) - 2(\bU - \sigma \mathbb{I}) \label{soluteL}
\end{align}
which holds on $I\times\Ozeta\subset \mathbb R^{1+2}$(with \eqref{shellEQL} posed on $I\times\omega\subset \mathbb R^{1+1}$)  where
\begin{align*}
\mathbb{S}_*= (\nabx \bv +(\nabx \bv)^\top) -\pi\mathbb{I}+  \bU,
\qquad
\mathbb{W}(\nabx \bv)=\tfrac{1}{2}(\nabx \bv-(\nabx \bv)^\top).
\end{align*}
We complement \eqref{divfreeL}--\eqref{soluteL} with the following initial conditions
\begin{align}
&\zeta(0,\cdot)=\zeta_0(\cdot), \qquad\partial_t\zeta(0,\cdot)=\zeta_\star(\cdot) & \text{in }\omega,
\\
&\bv(0,\cdot)=\bv_0(\cdot) & \text{in }\Omega_{\zeta_0},
\\
&\sigma(0,\cdot)=\sigma_0(\cdot),\quad\bU(0,\cdot)=\bU_0(\cdot) &\text{in }\Omega_{\zeta_0},
\label{initialCondSolvL}  
\end{align}
and again, we impose periodicity on the boundary of $\omega$ and the following interface condition
\begin{align} 
\label{interfaceL}
&\bv\circ\bm{\varphi}_\zeta=(\partial_t\zeta)\bn & \text{on }I\times \omega
\end{align}
at the flexible part of the boundary with normal $\bn$. We are interested in weak solutions of \eqref{divfreeL}--\eqref{interfaceL} where weak solutions are defined as follows:
\begin{definition}[Weak solution]
\label{def:weaksolmartFP}
Let $(\sigma_0, \bU_0, \bv_0, \zeta_0, \zeta_\star)$
be a dataset that satisfies
\begin{equation}
\begin{aligned}
\label{mainDataForAll}
%&\bff \in L^2(I;L^{2}_\mathrm{loc}(\mathbb{R}^2)),
%\qquad g\in L^2(I;L^{2}(\omega)) ,
%\\
&
\zeta_0 \in W^{2,2}(\omega) \text{ with } \Vert \zeta_0 \Vert_{L^\infty( \omega)} < L, \quad \zeta_\star \in L^{2}(\omega),
\\&\bv_0 \in L^{2}_{\divx}(\Omega_{\zeta_0} )\text{ is such that }\bv_0 \circ \bm{\varphi}_{\zeta_0} =\zeta_\star \bn \text{ on } \omega,
\\&
\sigma_0\in L^{2}(\Omega_{\zeta_0}), \quad
\bU_0\in L^{2}(\Omega_{\zeta_0}),
\\&
\sigma_0\geq 0,\,\, \bU_0>0 \quad \text{a.e. in } \Omega_{\zeta_0}.
\end{aligned}
\end{equation}
We call
$(\zeta, \bv, \sigma,\bU)$  
a \textit{weak solution} of   \eqref{divfreeL}--\eqref{interfaceL} with dataset $( \sigma_0, \bU_0, \bv_0, \zeta_0, \zeta_\star)$  if: 
\begin{itemize}
\item the following properties 
\begin{align*}
&\zeta\in  W^{1,\infty}\big(I;L^{2}(\omega)  \big) \cap W^{1,2}(I;W^{1,2}(\omega))  \cap L^{\infty}\big(I;W^{2,2}(\omega)  \big),
\\&
\Vert\zeta\Vert_{L^\infty(I\times\omega)}<L,
\\
&\bv\in
L^{\infty} \big(I; L^{2}(\Ozeta) \big)\cap L^2\big(I;W^{1,2}_{\divx}(\Ozeta)  \big),
\\
&
\sigma \in   L^{\infty}\big(I;L^{2}(\Ozeta)  \big),
\\
&
\bU \in   L^{\infty}\big(I;L^{2}(\Ozeta)  \big) ,
\\&
\sigma\geq 0,\,\, \bU>0 \quad \text{a.e. in } I\times\Ozeta
\end{align*}
holds;
\item  for all  $  \psi  \in C^\infty (\overline{I}\times \R^2  )$, we have
\begin{equation}
\begin{aligned} 
\label{weakRhoEq}
\int_I  \frac{\mathrm{d}}{\dt}
\int_{\Ozeta } \sigma\psi \dx \dt 
&=
\int_I
\int_{\Ozeta}[\sigma\partial_t\psi + (\sigma\bv\cdot\nabx) \psi] \dx\dt;
\end{aligned}
\end{equation}
\item for all  $  \mathbb{Y}  \in C^\infty (\overline{I}\times \R^2  )$, we have
\begin{equation}
\begin{aligned} 
\label{weakFokkerPlanckEq}
\int_I  \frac{\mathrm{d}}{\dt}
\int_{\Ozeta } \bU:\mathbb{Y} \dx \dt 
&=
\int_I
\int_{\Omega_{\zeta }}[\bU :\partial_t\mathbb{Y} + \bU:(\bv\cdot\nabx) \mathbb{Y}] \dx\dt
\\&+
\int_I\int_{\Ozeta}
[(\nabx \bv )\bU  + \bU (\nabx \bv )^\top]:\mathbb{Y} \dx\dt
\\&
-2\int_I\int_{\Ozeta }(\bU :\mathbb{Y}  - \sigma  \mathrm{tr}(\mathbb{Y} ))\dx\dt;
\end{aligned}
\end{equation}
\item for all  $(\bm{\phi},\phi)  \in C^\infty_{\divx} (\overline{I}\times \R^2  )\otimes  C^\infty (\overline{I}\times \omega )$ with $\bm{\phi}(T,\cdot)=0$, $\phi(T,\cdot)=0$ and $\bm{\phi}\circ \bm{\varphi}_\zeta= \phi \bn$, we have
\begin{equation}
\begin{aligned}
\label{weakFluidStrut}
\int_I  \frac{\mathrm{d}}{\dt}\bigg(
\int_{\Ozeta } \bv\cdot\bm{\phi} \dx 
+
\int_\omega\partial_t\zeta\phi\dy\bigg) \dt 
&=
\int_I
\int_{\Ozeta }[\bv \cdot\partial_t\bm{\phi} + \bv\cdot(\bv\cdot\nabx) \bm{\phi}] \dx\dt
\\&-
\int_I\int_{\Ozeta }\big[\nabx \bv  :\nabx \bm{\phi} +\bU:\nabx\bm{\phi}\big] \dx\dt
\\&
+
\int_I\int_\omega\big[\partial_t\zeta\partial_t\phi 
-
\gamma
\partial_t\naby\zeta\naby\phi 
-
\naby^2\zeta\naby^2\phi \big]\dy\dt;
\end{aligned}
\end{equation} 
%\item the energy inequality
%\begin{equation}
%\begin{aligned}   
%\sup_{t\in I}&
%\big(\Vert\sigma(t)\Vert_{L^2(\Ozeta)}^2
%+
%\Vert\bv(t)\Vert_{L^2(\Ozeta)}^2
%+
%\Vert\partial_t\zeta(t)\Vert_{L^2(\omega)}^2
%+
%\Vert\naby^2\zeta(t)\Vert_{L^2(\omega)}^2
%\big)
%\\
%&+
%\sup_{t\in I}
%\Vert\bU(t)\Vert_{L^2(\Ozeta)}^2 
%+
%\int_I\big(\Vert\nabx \bv \Vert_{L^2(\Ozeta)}^2
%+
%\gamma
%\Vert
%\partial_t\naby\zeta
%\Vert_{L^2(\omega)}^2
%+
%\Vert  \bU \Vert_{L^2(\Ozeta)}^2\big)\dt
%\\& 
%\lesssim
% (1+T)
%\Vert\sigma_0 \Vert_{L^2(\Omega_{\zeta_0})}^2
%+
%\Vert\bv_0 \Vert_{L^2(\Omega_{\zeta_0})}^2
%+
%\Vert \zeta_\star\Vert_{L^2(\omega)}^2
%+
%\Vert\naby^2\zeta_0\Vert_{L^2(\omega)}^2
%+
%\Vert\bU_0\Vert_{L^2(\Omega_{\zeta_0})}^2
%.
%\end{aligned}
%\end{equation} 
\end{itemize}
\end{definition}

%With this definition in hand, we now state our first main result.
%\begin{theorem}
%\label{thm:main}
%For a dataset $(\bff, g, \rho_0, \bT_0, \bu_0, \eta_0, \eta_\star)$
%that satisfies \eqref{mainDataForAll}, there  exists a global weak solution $(\eta, \bu, \rho,\bT)$ of   \eqref{divfree}--\eqref{interface}.  
%\end{theorem}

An apparent disadvantage of the system  \eqref{divfreeL}-\eqref{soluteL} without center-of-mass diffusion is that the regularisation effect due to said diffusion is lost. Consequently, the solute component $(\sigma,\bU)$ of its weak solution is less regular than the 
solute component $(\rho,\bT)$ of the weak solution of
\eqref{divfree}-\eqref{solute}, c.f. \cite{mensah2023weak}. Nevertheless, because  $(\sigma,\bU)$ possesses transport and corotational properties, any weak solution automatically regularises in the sense that it becomes essentially bounded in both space and time provided the initial conditions are equally bounded in space. Indeed, let $(\zeta, \bv, \sigma,\bU)$  
be a weak solution of   \eqref{divfreeL}--\eqref{interfaceL} with dataset $( \sigma_0, \bU_0, \bv_0, \zeta_0, \zeta_\star)$ satisfying \eqref{mainDataForAll} and $\sigma_0, \bU_0\in L^\infty(\Omega_{\zeta_0})$.
Equation \eqref{rhoEquL} being a transport equation means that it conserves all $L^p$-norms for $p\in[1,\infty)$, and thus, it is essentially bounded in spacetime. Although Equation \eqref{soluteL} looks complicated, we can show that it also dissipates all $L^p$-norms uniformly in $p\in [1,\infty)$ and thus, it is also essentially bounded in spacetime. To see the formal details, we test \eqref{rhoEquL} with $\sigma^{p-1}$, $p\in[1,\infty)$ and test \eqref{soluteL} with $\bU^{q-1}$, $q\in[1,\infty)$. For \eqref{rhoEquL}, we obtain for any $t\in I$,
\begin{align}\label{conserveSigma}
\Vert \sigma(t)\Vert_{L^p(\Omega_{\zeta})}^p=\Vert \sigma_0\Vert_{L^p(\Omega_{\zeta_0})}^p
\end{align}
uniformly in $p\in[1,\infty)$. Since $\sigma_0\in L^\infty(\Omega_{\zeta_0})$ by assumption, it follows that $\sigma \in L^\infty(I;L^\infty(\Omega_{\zeta}))$. 
\\
If we now test \eqref{soluteL} with $\bU^{q-1}$, $q\in[1,\infty)$ and use Proposition \ref{prop:zeroCorotational}, we obtain
\begin{align*}
\frac{1}{q}
\frac{\dd}{\dt}\Vert \bU\Vert_{L^q(\Ozeta)}^q
+2\Vert \bU\Vert_{L^q(\Ozeta)}^q
\leq
2\frac{1}{q}
\Vert \sigma\Vert_{L^q(\Ozeta)}^q
+
2
\frac{q-1}{q}
\Vert \bU\Vert_{L^q(\Ozeta)}^q .
\end{align*}
so that by \eqref{conserveSigma}, we obtain
\begin{align}\label{conserveSolute}
\Vert \bU(t)\Vert_{L^q(\Omega_{\zeta})}^q
+
2
\int_0^t
\Vert \bU\Vert_{L^q(\Omega_{\zeta})}^q
\ds
\leq
2T
\Vert \sigma_0\Vert_{L^q(\Omega_{\zeta_0})}^q 
+
\Vert \bU_0\Vert_{L^q(\Omega_{\zeta_0})}^q 
\end{align}
uniformly in $q\in[1,\infty)$. So, in particular, we obtain $\bU \in L^\infty(I;L^\infty(\Omega_{\zeta}))$ since $\sigma_0,\bU_0\in L^\infty(\Omega_{\zeta_0})$ by assumption. 

The discussion above motivates a finer notion (for the solute component) of a weak solution of   \eqref{divfreeL}--\eqref{interfaceL}. We shall refer to this as an \textit{essentially bounded weak solution} whose precise definition is given as follows:
\begin{definition}[Essentially bounded weak solution]
\label{def:weaksolmartFPBounded}
Let $(\sigma_0, \bU_0, \bv_0, \zeta_0, \zeta_\star)$
be a dataset that satisfies \eqref{mainDataForAll} and
\begin{equation}
\begin{aligned}
\label{mainDataForAllBounded}
\sigma_0\in L^{\infty}(\Omega_{\zeta_0}), \quad
\bU_0\in L^{\infty}(\Omega_{\zeta_0}).
\end{aligned}
\end{equation}
We call
$(\zeta, \bv, \sigma,\bU)$  
an \textit{essentially bounded weak solution} of   \eqref{divfreeL}--\eqref{interfaceL} with dataset $( \sigma_0, \bU_0, \bv_0, \zeta_0, \zeta_\star)$  if:
\begin{itemize}
\item $(\zeta, \bv, \sigma,\bU)$  
is a  weak solution of   \eqref{divfreeL}--\eqref{interfaceL} with dataset $( \sigma_0, \bU_0, \bv_0, \zeta_0, \zeta_\star)$;
\item the pair $(\sigma,\bU)$ satisfies $
\sigma, \bU \in L^\infty(I;L^\infty(\Omega_{\zeta})).
$
\end{itemize}
\end{definition}
With Definition \ref{def:weaksolmartFPBounded} in hand, we now state our first main theorem on the existence of an essentially bounded weak solution of   \eqref{divfreeL}--\eqref{interfaceL}.
\begin{theorem}
\label{thm:main} 
For a dataset $(\sigma_0, \bU_0, \bv_0, \zeta_0, \zeta_\star)$ satisfying \eqref{mainDataForAll}
with $(\sigma_0,\bU_0)$ satisfying \eqref{mainDataForAllBounded}, there exists an essentially bounded weak solution $(\zeta,\bv,  \sigma,\bU)$  of \eqref{divfreeL}-\eqref{interfaceL}.
\end{theorem}
Our second main result will now seek to establish a relationship between strong solutions of \eqref{divfree}-\eqref{interface} and essentially bounded weak solutions of   \eqref{divfreeL}--\eqref{interfaceL}.  
This vanishing center-of-mass limit result was earlier announced in \cite{masmoudi2013global} for the system posed on a fixed domain and for the mesoscopic description \eqref{fokkerPlanck} of the solute. See the first bullet point on \cite[Page 486]{masmoudi2013global}.
This result was, however, not released by the author. This open problem is revisited in \cite[Page 336]{debiec2023corotational} with a brief discussion on \cite[Section 6.4]{debiec2023corotational}. The system the authors discuss in the latter result, however, has an additional symmetric gradient term $\tfrac{b}{2}(\nabx\bu+(\nabx\bu)^\top)$, $b\geq0$ in the equation for the extra stress tensor $\bT$ and as such, does not arise as the macroscopic closure of \eqref{fokkerPlanck} when $b>0$. Nevertheless, from the purely mathematical point of view, this additional term provides a cancellation property in the energy identity for the extra stress term $\divx\bT$ in the momentum equation \eqref{momEq}.

In this work, we consider the actual macroscopic closure of \eqref{fokkerPlanck} and we do not require any assumption of the tensor gradient $\nabx\bU$ as discussed in \cite[Section 6.4]{debiec2023corotational}. The precise statement is given as follows.
\begin{theorem}
\label{thm:main2}
Let $(\sigma_0, \bU_0, \bv_0, \zeta_0, \zeta_\star)$ be a dataset satisfying \eqref{mainDataForAll} and \eqref{mainDataForAllBounded} and let $({\rho}^\varepsilon_0, {\bT}^\varepsilon_0, {\bu}^\varepsilon_0, \eta^\varepsilon_0, \eta^\varepsilon_\star)_{\varepsilon>0}$ be datasets satisfying \eqref{mainDataForAllStrong} and
\begin{equation}
\begin{aligned}
\label{dataConv}
&\eta^\varepsilon_0 \rightarrow \zeta_0 \quad\text{in}\quad W^{2,2}( \omega), 
\\&\eta^\varepsilon_\star \rightarrow \zeta_\star \quad\text{in}\quad L^\infty(I;L^{2}( \omega)), 
\\&\bm{1}_{\Omega_{\eta^\varepsilon_0}} \bu^\varepsilon_0  \rightarrow \bm{1}_{\Omega_{\zeta_0}}\bv_0 \quad\text{in}\quad L^2(\Omega \cup S_\ell),  
\\&\bm{1}_{\Omega_{\eta^\varepsilon_0}} \rho^\varepsilon_0 \rightarrow \bm{1}_{\Omega_{\zeta_0}}\sigma_0 \quad\text{in}\quad L^2(\Omega \cup S_\ell)
, 
\\&\bm{1}_{\Omega_{\eta^\varepsilon_0}}
\bT^\varepsilon_0  \rightarrow \bm{1}_{\Omega_{\zeta_0}}\bU_0 \quad\text{in}\quad L^2(\Omega \cup S_\ell)
\end{aligned}
\end{equation}
as $\varepsilon\rightarrow0$.  
Now let $(\eta^\varepsilon, {\bu}^\varepsilon, {p}^\varepsilon, {\rho}^\varepsilon, {\bT}^\varepsilon)_{\varepsilon>0}$ be a collection of strong solutions of \eqref{divfree}-\eqref{interface} with dataset $({\rho}^\varepsilon_0, {\bT}^\varepsilon_0, {\bu}^\varepsilon_0, \eta^\varepsilon_0, \eta^\varepsilon_\star)_{\varepsilon>0}$ and let $(\zeta,\bv,  \sigma,\bU)$ be  an essentially bounded weak solution of \eqref{divfreeL}-\eqref{interfaceL} with dataset $(\sigma_0, \bU_0, \bv_0, \zeta_0, \zeta_\star)$. Then the following convergences
\begin{align*}
&\eta^\varepsilon \rightarrow \zeta \quad\text{in}\quad L^\infty(I;W^{2,2}( \omega)), 
\\&\partial_t\eta^\varepsilon \rightarrow \partial_t\zeta \quad\text{in}\quad L^\infty(I;L^{2}( \omega)),
\\&\partial_t\eta^\varepsilon \rightarrow \partial_t\zeta \quad\text{in}\quad L^2(I;W^{1,2}( \omega)),
\\&\bm{1}_{\Omega_{\eta^\varepsilon}} \bu^\varepsilon  \rightarrow \bm{1}_{\Omega_{\zeta}}\bv \quad\text{in}\quad L^\infty(I;L^2(\Omega \cup S_\ell)), 
\\&\bm{1}_{\Omega_{\eta^\varepsilon}}\nabx\bu^\varepsilon  \rightarrow \bm{1}_{\Omega_{\zeta}}\nabx\bv \quad\text{in}\quad L^2(I;L^2(\Omega \cup S_\ell)), 
\\&\bm{1}_{\Omega_{\eta^\varepsilon}} \rho^\varepsilon \rightarrow \bm{1}_{\Omega_{\zeta}}\sigma \quad\text{in}\quad L^\infty(I;L^2(\Omega \cup S_\ell))
, 
\\&\bm{1}_{\Omega_{\eta^\varepsilon}}
\bT^\varepsilon  \rightarrow \bm{1}_{\Omega_{\zeta}}\bU \quad\text{in}\quad L^\infty(I;L^2(\Omega \cup S_\ell))
\end{align*}
hold. 
\end{theorem}
\begin{remark}
Theorem \ref{thm:main2} trivially holds for fixed spatial domains where formally speaking, $\eta^\varepsilon=\zeta=0$ and $\Omega_{\eta^\varepsilon}=\Ozeta=\Omega$.
\end{remark} 
An immediate consequence of Theorem \ref{thm:main2} is a \textit{weak-strong uniqueness} results. Generally speaking, such a result states that  a weak solution is unique in the class of a strong solution.
This will usually require that both solutions solves the same equation and may also require additionally regularity of either solution or their data. See for instant \cite{BMSS, chemetov2019weak, germain2011weak, schwarzacher2022weak}. 
In our setting, however, \eqref{divfree}--\eqref{interface} and \eqref{divfreeL}--\eqref{interfaceL}  are two different systems, albeit related, with their own notion of a solution. Nevertheless, we present a \textit{`limiting' weak-strong uniqueness} result which states that as $\varepsilon\rightarrow0$, the essentially bounded weak solution to \eqref{divfreeL}--\eqref{interfaceL} will be unique in the class of the strong solution to \eqref{divfree}--\eqref{interface}. Here, we do not require additional regularity assumption on their respective dataset. The precise statement is:
 \begin{proposition}
Let $(\eta,\bu, p,\rho,\bT)$  be a strong solution of \eqref{divfree}-\eqref{interface} with dataset $(\rho_0, \bT_0, \bu_0, \eta_0, \eta_\star)$ satisfying \eqref{mainDataForAllStrong} and let $(\zeta,\bv,  \sigma,\bU)$ be an essentially bounded weak solution of \eqref{divfreeL}-\eqref{interfaceL} with dataset $(\sigma_0, \bU_0, \bv_0, \zeta_0, \zeta_\star)$ satisfying \eqref{mainDataForAll} and \eqref{mainDataForAllBounded}.
Define
\begin{align*}
\overline{\bu} =\bu \circ \bm{\Psi}_{\eta-\zeta}, \quad \overline{p} =p \circ \bm{\Psi}_{\eta-\zeta},\quad \overline{\rho} =\rho \circ \bm{\Psi}_{\eta-\zeta}, \quad \overline{\bT} =\bT \circ \bm{\Psi}_{\eta-\zeta}
\end{align*}
with respect to the  Hanzawa transform $\bm{\Psi}_{\eta-\zeta}$.
% If either
%\begin{enumerate}
%\item $\varepsilon \rightarrow0$ or;
%\item $\overline{\rho}\neq\sigma$ and $\frac{\Vert (\overline{\rho}-\sigma)(t)\Vert_{L^\infty(\Omega_{\zeta})}}{
%\Vert (\overline{\rho}-\sigma)(t)\Vert_{L^2(\Omega_{\zeta})} 
%}=\mathcal{O}(1)$ for a.e. $t>0$;\footnote{$\mathcal{O}(1)$ is the Landau notations for quantities that are bounded}
%\end{enumerate}
%then
As $\varepsilon \rightarrow0$, we have
\begin{equation}
\begin{aligned} 
\sup_{t'\in (0,t]}&\big(
 \Vert  \partial_t(\eta-\zeta)(t') \Vert_{L^{2}(\omega )}^2
+
\Vert  \naby^2(\eta-\zeta)(t') \Vert_{L^{2}(\omega )}^2 
+
\Vert (\overline{\bu}-\bv)(t')\Vert_{L^2(\Omega_{\zeta})}^2
\big)
\\&
+\sup_{t'\in (0,t]}\big(\Vert (\overline{\rho}-\sigma)(t')\Vert_{L^2(\Omega_{\zeta})}^2
+
\Vert (\overline{\bT}-\bU)(t')\Vert_{L^2(\Omega_{\zeta})}^2
\big)
+
\int_0^t
\Vert \overline{\bT}-\bU \Vert_{L^2(\Omega_{\zeta})}^2\dt'
\\&+ 
\int_0^t
\Vert \nabx(\overline{\bu}-\bv)\Vert_{L^2(\Omega_{\zeta})}^2\dt'
+ 
\gamma 
\int_0^t\Vert  \partial_{t'}\naby(\eta-\zeta) \Vert_{L^{2}(\omega )}^2
\dt'  
\\
\lesssim& 
 \Vert   \eta_\star-\zeta_\star \Vert_{L^{2}(\omega )}^2
+
\Vert  \naby^2(\eta_0-\zeta_0) \Vert_{L^{2}(\omega )}^2  
+
\Vert \overline{\bu}_0-\bv_0\Vert_{L^2(\Omega_{\zeta_0})}^2
\\&
\qquad+
 \Vert \overline{\rho}_0-\sigma_0\Vert_{L^2(\Omega_{\zeta_0})}^2  
+
\Vert \overline{\bT}_0-\bU_{0}\Vert_{L^2(\Omega_{\zeta_0})}^2 
\end{aligned}
\end{equation}
for all $t\in I$ with a constant that depends on  $\gamma$ and $T$.
\end{proposition}
\begin{proof}
As mentioned earlier, the proof of this result is a direct consequence of Theorem \ref{thm:main2} and in particular, achieved by passing to the limit $\varepsilon\rightarrow0$ in \eqref{contrEst0} of Proposition \ref{prop:main} below.
\end{proof}
\begin{remark}
Unfortunately we are unable to recover actual unconditional uniqueness for essentially bounded weak solutions of \eqref{divfreeL}-\eqref{interfaceL} with dataset   satisfying \eqref{mainDataForAll} and \eqref{mainDataForAllBounded} from our analysis and it also not clear if that is to be expected.  
\end{remark}

\subsection{Plan for the rest of the paper}
In the next section, Section \ref{sec:prelim}, we collect some notations, set up our functional framework, and collect some key results that would be used repeatedly in the proof of our main results. We then devote Section \ref{sec:Existence} and Section \ref{sec:singularLimit} to the proofs of Theorem \ref{thm:main} and Theorem \ref{thm:main2}, respectively.

\section{Preliminaries }
\label{sec:prelim}
\noindent  
%Henceforth, without loss of generality, we set all the parameters $\{ \varepsilon,K,\gamma,k,\nu,\varrho_s \}$ in \eqref{divfree}-\eqref{bddCondSolv} to 1.
For two non-negative quantities $F$ and $G$, we write $F \lesssim G$  if there is a constant $c>0$  such that $F \leq c\,G$. If $F \lesssim G$ and $G\lesssim F$ both hold, we use the notation $F\sim G$.  The anti-symmetric gradient of a vector $\bff\in \mathbb{R}^d$ is denoted by $\mathbb{W}(\nabx \bff)=\tfrac{1}{2}(\nabx \bff-(\nabx \bff)^\top)$ and the scalar matrix product of the matrices $\mathbb{A}=(a_{ij})_{i,j=1}^d$ and $\mathbb{B}=(b_{ij})_{i,j=1}^d$ is denoted by $\mathbb{A}:\mathbb{B}=\sum_{ij}a_{ij}b_{ji}$.
The symbol $\vert \cdot \vert$ may be used in four different contexts. For a scalar function $f\in \mathbb{R}$, $\vert f\vert$ denotes the absolute value of $f$. For a vector $\bff\in \mathbb{R}^d$, $\vert \bff \vert$ denotes the Euclidean norm of $\bff$. For a square matrix $\mathbb{F}\in \mathbb{R}^{d\times d}$, $\vert \mathbb{F} \vert$ shall denote the Frobenius norm $\sqrt{\mathrm{trace}(\mathbb{F}^T\mathbb{F})}$. Also, if $S\subseteq  \mathbb{R}^d$ is  a (sub)set, then $\vert S \vert$ is the $d$-dimensional Lebesgue measure of $S$.  
Since we only consider functions on $\omega \subset\mathbb{R}$ with periodic boundary
conditions and zero mean values, we have the following equivalences
\begin{align}
\label{equiNorm}
\Vert \cdot\Vert_{W^{1,2}(\omega)}\sim
\Vert \partial_y\cdot\Vert_{L^{2}(\omega)},
\qquad
\Vert \cdot\Vert_{W^{2,2}(\omega)}\sim
\Vert \partial_y^2\cdot\Vert_{L^{2}(\omega)},
\qquad
\Vert \cdot\Vert_{W^{4,2}(\omega)}\sim
\Vert \partial_y^4\cdot\Vert_{L^{2}(\omega)}.
\end{align} 
For $I:=(0,T)$, $T>0$, and $\eta\in C(\overline{I}\times\omega)$ satisfying $\|\eta\|_{L^\infty(I\times\omega)}\leq L$ where $L>0$ is a constant, we define for $1\leq p,r\leq\infty$,
\begin{align*} 
L^p(I;L^r(\Omega_\eta))&:=\Big\{v\in L^1(I\times\Omega_\eta):\substack{v(t,\cdot)\in L^r(\Omega_{\eta(t)})\,\,\text{for a.e. }t,\\\|v(t,\cdot)\|_{L^r(\Omega_{\eta(t)})}\in L^p(I)}\Big\},\\
L^p(I;W^{1,r}(\Omega_\eta))&:=\big\{v\in L^p(I;L^r(\Omega_\eta)):\,\,\nabx v\in L^p(I;L^r(\Omega_\eta))\big\}.
\end{align*} 
Higher-order Sobolev spaces can be defined accordingly. For $k>0$ with $k\notin\mathbb N$, we define the fractional Sobolev space $L^p(I;W^{k,r}(\Oeta))$ as the class of $L^p(I;L^r(\Omega_\eta))$-functions $v$ for which 
\begin{align*}
\|v\|_{L^p(I;W^{k,r}(\Oeta))}^p
&=\int_I\bigg(\int_{\Oeta} \vert v\vert^r\dx
+\int_{\Oeta}\int_{\Oeta}\frac{|v(\bx)-v(\bx')|^r}{|\bx-\bx'|^{d+k r}}\dx\dx'\bigg)^{\frac{p}{r}}\dt
\end{align*}
is finite. Accordingly, we can also introduce fractional differentiability in time for the spaces on moving domains.

Next, for  $\eta\in C(\overline{I}\times\omega)$ satisfying $\|\eta\|_{L^\infty(I \times \omega)}\leq L$ and  $\|\naby\eta\|_{L^\infty(I \times \omega)}\leq c(L)$,  where $L>0$ is a constant, we let
\begin{align*}
\mathrm{Bog}_\eta : C_0^\infty(\Oeta) \rightarrow C_0^\infty(\Oeta; \mathbb{R}^3) \quad\text{with}\quad \divx \mathrm{Bog}_\eta(f)= f-b\int_{\Oeta}f\dx,
\end{align*}
with $b\in C_0^\infty(\Oeta\setminus S_L)$, be the time-dependent Bogovskij operator (see \cite[Theorem 2.11 \& Remark 2.12]{BMSS}). The operator $\mathrm{Bog}_\eta$ vanishes on the boundary $\partial\Oeta$, it commutes with time derivatives, and it continuously maps scalar elements in $W^{s,p}(\Ozeta)$ into vectors in  $W^{s+1,p}(\Ozeta;\mathbb{R}^3)$ for any $p\in(1,\infty)$ and $s\geq0$. 
Furthermore,  the Hanzawa transform $\bm{\Psi}_\eta$ together with its inverse 
 $\bm{\Psi}_\eta^{-1} : \Oeta \rightarrow\Omega$  possesses the following properties, see \cite{breit2022regularity, BMSS} for details. If for some $\ell,R>0$, we assume that
\begin{align*}
\Vert\eta\Vert_{L^\infty(\omega)}
+
\Vert\zeta\Vert_{L^\infty(\omega)}
< \ell <L \qquad\text{and}\qquad
\Vert\naby\eta\Vert_{L^\infty(\omega)}
+
\Vert\naby\zeta\Vert_{L^\infty(\omega)}
<R
\end{align*}
holds, then for any  $s>0$, $\varrho,p\in[1,\infty]$ and for any $\eta,\zeta \in B^{s}_{\varrho,p}(\omega)\cap W^{1,\infty}(\omega)$ (where $B^{s}_{\varrho,p}$ is a Besov space), we have that the estimates
\begin{align}
\label{210and212}
&\Vert \bm{\Psi}_\eta \Vert_{B^s_{\varrho,p}(\Omega\cup S_\ell)}
+
\Vert \bm{\Psi}_\eta^{-1} \Vert_{B^s_{\varrho,p}(\Omega\cup S_\ell)}
 \lesssim
1+ \Vert \eta \Vert_{B^s_{\varrho,p}(\omega)},
\\
\label{211and213}
&\Vert \bm{\Psi}_\eta - \bm{\Psi}_\zeta  \Vert_{B^s_{\varrho,p}(\Omega\cup S_\ell)} 
+
\Vert \bm{\Psi}_\eta^{-1} - \bm{\Psi}_\zeta^{-1}  \Vert_{B^s_{\varrho,p}(\Omega\cup S_\ell)} 
\lesssim
 \Vert \eta - \zeta \Vert_{B^s_{\varrho,p}(\omega)}
\end{align}
and
\begin{align}
\label{218}
&\Vert \partial_t\bm{\Psi}_\eta \Vert_{B^s_{\varrho,p}(\Omega\cup S_\ell)}
\lesssim
 \Vert \partial_t\eta \Vert_{B^{s}_{ \varrho,p}(\omega)},
\qquad
\eta
\in W^{1,1}(I;B^{s}_{\varrho,p}(\omega))
\end{align}
holds uniformly in time with the hidden constants depending only on the reference geometry, on $L-\ell$ and $R$. Finally, before we prove our main result, we present a useful result we shall use at several points in our proof.
\begin{proposition}
\label{prop:zeroCorotational}
Let $\bw=(w_1,w_2)$ be a $2$-d vector and $\mathbb{Z}=(z_{ij})_{i,j=1}^2$ a $2\times 2$ matrix. Then for any $n\in\mathbb{N}$ and for any $\mathbb{Y}\in\{\mathbb{Z},\mathbb{Z}^\top\}$, the equation
\begin{align*}
\mathbb{W}(\nabx \bw)\mathbb{Z}:\mathbb{Y}
^n + \mathbb{Z}\mathbb{W}((\nabx \bw)^\top) :\mathbb{Y}^n=0   
\end{align*}
holds where $\mathbb{Y}^n$ is the $n$-times matrix multiplication.\footnote{Note that the same result holds for  $\mathbb{A}:\mathbb{B}:=\sum a_{ij}b_{ij}$ rather than $\sum a_{ij}b_{ji}$ as stated in Section \ref{sec:prelim}.}
%, $\mathbb{W}(\nabx \bw)=\tfrac{1}{2}(\nabx \bw-(\nabx \bw)^\top)$ is the anti-symmetric gradient of $\bw$ and $\mathbb{A}:\mathbb{B}=\sum_{ij}a_{ij}b_{ji}$ is the scaler matrix product of the matrices $\mathbb{A}=(a_{ij})_{i,j=1}^2$ and $\mathbb{B}=(b_{ij})_{i,j=1}^2$.
\end{proposition}
\begin{proof}
We consider the case $\mathbb{Y}=\mathbb{Z}$ with $\mathbb{Y}=\mathbb{Z}^\top$ following a similar process. 
\\
The desired equation can be proved by induction and so we first take $n=1$. A straightforward calculation gives
\begin{align*}
2\mathbb{W}(\nabx \bw)\mathbb{Z} 
&=
\begin{pmatrix}
0 &\partial_2w_1-\partial_1w_2\\
-(\partial_2w_1-\partial_1w_2) & 0
\end{pmatrix}
\begin{pmatrix}
z_{11}&z_{12}\\
z_{21}&z_{22}
\end{pmatrix}
\\
&=
\begin{pmatrix}
(\partial_2w_1-\partial_1w_2)z_{21} &(\partial_2w_1-\partial_1w_2)z_{22}\\
-(\partial_2w_1-\partial_1w_2)z_{11} & -(\partial_2w_1-\partial_1w_2)z_{12}
\end{pmatrix}
\end{align*}
therefore,
\begin{align*}
2\mathbb{W}(\nabx \bw)\mathbb{Z}:\mathbb{Z}
=
(\partial_2w_1-\partial_1w_2)(z_{11}z_{21}+z_{21}z_{22}-z_{11}z_{12}-z_{12}z_{22}).
\end{align*}
Also,
\begin{align*}
2\mathbb{Z} \mathbb{W}((\nabx \bw)^\top)
&=
\begin{pmatrix}
z_{11}&z_{12}\\
z_{21}&z_{22}
\end{pmatrix}
\begin{pmatrix}
0 &-(\partial_2w_1-\partial_1w_2)\\
\partial_2w_1-\partial_1w_2 & 0
\end{pmatrix}
\\
&=
\begin{pmatrix}
(\partial_2w_1-\partial_1w_2)z_{12} &-(\partial_2w_1-\partial_1w_2)z_{11}\\
(\partial_2w_1-\partial_1w_2)z_{22} & -(\partial_2w_1-\partial_1w_2)z_{21}
\end{pmatrix}
\end{align*}
and therefore,
\begin{align*}
2\mathbb{Z} \mathbb{W}((\nabx \bw)^\top):\mathbb{Z}
=&-
(\partial_2w_1-\partial_1w_2)(z_{11}z_{21}+z_{21}z_{22}-z_{11}z_{12}-z_{12}z_{22})
\\=&-2\mathbb{W}(\nabx \bw)\mathbb{Z}:\mathbb{Z}.
\end{align*}
If we now assume that $n=k$ true, then for $n=k+1$, we obtain
\begin{align*}
\mathbb{W}(\nabx \bw)\mathbb{Z}:\mathbb{Z}
^{k+1}+ \mathbb{Z}\mathbb{W}((\nabx \bw)^\top) :\mathbb{Z}^{k+1}
=
-\mathbb{Z}\mathbb{W}((\nabx \bw)^\top) \mathbb{Z}^{k}\mathbb{Z}+ \mathbb{Z}\mathbb{W}((\nabx \bw)^\top) :\mathbb{Z}^{k}\mathbb{Z}=0.
\end{align*}
This finishes the proof.
\end{proof}

\section{Existence of essentially bounded weak solutions}
\label{sec:Existence}
We devote this section to the proof of our second main result, Theorem \ref{thm:main}. To begin with, we consider a diffusive regularisation of \eqref{divfreeL}-\eqref{interfaceL} with  Laplacian terms on the right of \eqref{rhoEquL} and \eqref{soluteL} so that equations for the solute subproblem is no different from the solute subproblem of  \eqref{divfree}-\eqref{interface}. More precisely, we consider the following problem
\begin{align}
\divx \bv=0, \label{divfreeLReg}
\\
\partial_t \sigma+ (\bv\cdot \nabx) \sigma
= 
\alpha\Delx\sigma
,\label{rhoEquLReg}
\\
\partial_t \bv  + (\bv\cdot \nabx)\bv
= 
  \delx \bv -\nabx \pi 
+ 
\divx   \bU, \label{momEqLReg}
\\
 \partial_t^2 \zeta  - \gamma\partial_t\partial_y^2 \zeta +\partial_y^4 \zeta =  - ( \mathbb{S}_*\bn_\zeta )\circ \bm{\varphi}_\zeta\cdot\bn \,\det(\partial_y\bm{\varphi}_\zeta), 
\label{shellEQLReg}
\\
\partial_t \bU + (\bv\cdot \nabx) \bU
=
\mathbb{W}(\nabx \bv)\bU + \bU\mathbb{W}((\nabx \bv)^\top) - 2(\bU - \sigma \mathbb{I})
+\alpha\Delx\bU \label{soluteLReg}
\end{align}
posed on $I\times\Ozeta\subset \mathbb R^{1+2}$ where
\begin{align*}
\mathbb{S}_*= (\nabx \bv +(\nabx \bv)^\top) -\pi\mathbb{I}+  \bU,
\qquad
\mathbb{W}(\nabx \bv)=\tfrac{1}{2}(\nabx \bv-(\nabx \bv)^\top).
\end{align*}
We complement \eqref{divfreeLReg}--\eqref{soluteLReg} with the following initial and boundary conditions
\begin{align}
&\zeta(0,\cdot)=\zeta_0(\cdot), \qquad\partial_t\zeta(0,\cdot)=\zeta_\star(\cdot) & \text{in }\omega,
\\
&\bv(0,\cdot)=\bv_0(\cdot) & \text{in }\Omega_{\zeta_0},
\\
&\sigma(0,\cdot)=\sigma_0(\cdot),\quad\bU(0,\cdot)=\bU_0(\cdot) &\text{in }\Omega_{\zeta_0},
\label{initialCondSolvLReg}  
\\
& 
\bn_{\zeta}\cdot\nabx\sigma=0,\qquad
\bn_{\zeta}\cdot\nabx\bU=0 &\text{on }I\times\partial\Omega_{\zeta}.
\label{bddCondSolvZeta}
\end{align}
and again, we impose periodicity on the boundary of $\omega$ and the following interface condition
\begin{align} 
\label{interfaceLReg}
&\bv\circ\bm{\varphi}_\zeta=(\partial_t\zeta)\bn & \text{on }I\times \omega
\end{align}
at the flexible part of the boundary with normal $\bn$.

A \textit{weak solution} $(\zeta, \bv, \sigma,\bU)$  
 of   \eqref{divfreeLReg}--\eqref{interfaceLReg}, in the sense of Definition \ref{def:weaksolmartFP}, emanating from a dataset $( \sigma_0, \bU_0, \bv_0, \zeta_0, \zeta_\star)$ satisfying \eqref{mainDataForAll} has been proved in  \cite{mensah2023weak}. The system  \eqref{divfreeLReg}--\eqref{interfaceLReg} being corotational (i.e. having the vorticity gradient $\mathbb{W}(\nabx \bv)$ on the right of \eqref{soluteLReg} rather than the full gradient $\nabx \bv$) actually simplifies the proof due to Proposition \ref{prop:zeroCorotational}.
For completeness, but to avoid repetition, we give a formal proof of the construction of a \textit{weak solution} (see \cite[Theorem 2.2]{mensah2023weak} for the details) as well as how such a weak solution can be extended to an \textit{essentially bounded weak solution} given that $(\sigma_0,\bU_0)$ additionally satisfies \eqref{mainDataForAllBounded}. A rigorous justification will then follow by performing the subsequent analysis on the Galerkin level and passing to the limit. 

To obtain a weak solution, we consider $(\sigma,\bv, \partial_t\zeta,\bU)$ as test functions for \eqref{rhoEquLReg}-\eqref{soluteLReg}, respectively. This yields
\begin{align*}
&\frac{\dd}{\dt}\Vert\sigma\Vert_{L^2(\Ozeta)}^2+\alpha\Vert\nabx\sigma\Vert_{L^2(\Ozeta)}^2=0,
\\
&\frac{\dd}{\dt}\big(\Vert\bv\Vert_{L^2(\Ozeta)}^2 + \Vert\partial_t\zeta\Vert_{L^2(\omega)}^2
+ \Vert\naby^2\zeta\Vert_{L^2(\omega)}^2
\big)+\Vert\nabx\bv\Vert_{L^2(\Ozeta)}^2
+
\gamma\Vert \partial_t\naby\zeta\Vert_{L^2(\omega)}^2 =0,
\\
&\frac{\dd}{\dt}\Vert\bU\Vert_{L^2(\Ozeta)}^2+\alpha\Vert\nabx\bU\Vert_{L^2(\Ozeta)}^2 + \Vert\bU\Vert_{L^2(\Ozeta)}^2
\leq \Vert\sigma\Vert_{L^2(\Ozeta)}^2,
\end{align*}
where we use Proposition \ref{prop:zeroCorotational} to deal with the first two terms on the right-hand side of \eqref{soluteLReg}. Consequently, the requirement of estimating the trace of \eqref{soluteLReg} as was done in \cite[Proposition 3.1]{mensah2023weak}  is no longer needed. We now obtain by integrating the three (in)equalities above in time,
\begin{equation}
\begin{aligned}   
\label{3.11}
\sup_{t\in I}&
\big(\Vert\sigma(t)\Vert_{L^2(\Ozeta)}^2
+
\Vert\bv(t)\Vert_{L^2(\Ozeta)}^2
+
\Vert\partial_t\zeta(t)\Vert_{L^2(\omega)}^2
+
\Vert\naby^2\zeta(t)\Vert_{L^2(\omega)}^2
\big)
\\
&+
\sup_{t\in I}
\Vert\bU(t)\Vert_{L^2(\Ozeta)}^2 
+
\int_I\big(\Vert\nabx \bv \Vert_{L^2(\Ozeta)}^2
+
\gamma\Vert \partial_t\naby\zeta\Vert_{L^2(\omega)}^2
\big)\dt
\\&
+
\int_I\big(
\alpha\Vert\nabx\sigma\Vert_{L^2(\Ozeta)}^2
+\alpha\Vert\nabx\bU\Vert_{L^2(\Ozeta)}^2 
+
\Vert  \bU \Vert_{L^2(\Ozeta)}^2\big)\dt
\\& 
\lesssim
(1+T)
\Vert\sigma_0 \Vert_{L^2(\Omega_{\zeta_0})}^2
+
\Vert\bv_0 \Vert_{L^2(\Omega_{\zeta_0})}^2
+
\Vert \zeta_\star\Vert_{L^2(\omega)}^2
+
\Vert\naby^2\zeta_0\Vert_{L^2(\omega)}^2
+
\Vert\bU_0\Vert_{L^2(\Omega_{\zeta_0})}^2
.
\end{aligned}
\end{equation} 
with a constant that is independent of $\alpha>0$. Consequently, we can pass to the limit $\alpha \rightarrow0$ and obtain a weak solution $(\zeta, \bv, \sigma,\bU)$   of   \eqref{divfreeL}--\eqref{interfaceL}. The fact that $(\zeta, \bv, \sigma,\bU)$ is an essentially bounded weak solution  of   \eqref{divfreeL}--\eqref{interfaceL} when $(\sigma_0,\bU_0)$ also satisfies \eqref{mainDataForAllBounded} has already been explained in \eqref{conserveSigma} and \eqref{conserveSolute}. This finishes the proof. 
\begin{remark}
We need to clarify that the formal passage to the limit $\alpha\rightarrow0$ in \eqref{3.11} for the construction of essentially bounded weak solutions is different from Theorem \ref{thm:main2} whose proof is to be presented in Section \ref{sec:singularLimit} below. The former is just an approximation layer on top of the finite-dimensional Galerkin approximation that allows for compactness. The latter is a singular limit result that compares two different systems with their respective notion of a solution and shows that a family of solutions to one system converges to the solution of the other.
\end{remark} 

\section{Vanishing center-of-mass limit}
\label{sec:singularLimit}
We devote this section to the proof of our second main result, Theorem \ref{thm:main2}.
Our proof relies on a relative energy method where one measures the distance between two solutions in the energy norm. Unfortunately, because our two solutions under consideration are defined on two separate variable geometries, we are unable to directly measure their distance. To get around this issue, we transform one solution via the  Hanzawa transform onto the domain of the other and measure the distance on this latter domain. In this regard, we consider the following transformation 
\begin{align}
\label{transHan}
\overline{\bu} =\bu \circ \bm{\Psi}_{\eta-\zeta}, \quad \overline{p} =p \circ \bm{\Psi}_{\eta-\zeta},\quad \overline{\rho} =\rho \circ \bm{\Psi}_{\eta-\zeta}, \quad \overline{\bT} =\bT \circ \bm{\Psi}_{\eta-\zeta}
\end{align}
of  $(\eta ,\bu , p ,\rho ,\bT )$ onto the domain of  $(\zeta, \bv, \pi, \sigma,\bU)$ with respect to the Hanzawa transform $\bm{\Psi}_{\eta-\zeta}$. The proof of Theorem \ref{thm:main2} will now follow from the following proposition.
\begin{proposition}\label{prop:main}
Let $(\eta,\bu, p,\rho,\bT)$  be a strong solution of \eqref{divfree}-\eqref{interface} with dataset $(\rho_0, \bT_0, \bu_0, \eta_0, \eta_\star)$. If $(\zeta,\bv,  \sigma,\bU)$ is an essentially bounded weak solution of \eqref{divfreeL}-\eqref{interfaceL} with dataset $(\sigma_0, \bU_0, \bv_0, \zeta_0, \zeta_\star)$, then we have
\begin{equation}
\begin{aligned}
\label{contrEst0} 
\sup_{t'\in (0,t]}&\big(
 \Vert  \partial_t(\eta-\zeta)(t') \Vert_{L^{2}(\omega )}^2
+
\Vert  \naby^2(\eta-\zeta)(t') \Vert_{L^{2}(\omega )}^2 
+
\Vert (\overline{\bu}-\bv)(t')\Vert_{L^2(\Omega_{\zeta})}^2
\big)
\\&
+\sup_{t'\in (0,t]}\big(\Vert (\overline{\rho}-\sigma)(t')\Vert_{L^2(\Omega_{\zeta})}^2
+
\Vert (\overline{\bT}-\bU)(t')\Vert_{L^2(\Omega_{\zeta})}^2
\big)
+
\int_0^t
\Vert \overline{\bT}-\bU \Vert_{L^2(\Omega_{\zeta})}^2\dt'
\\&+ 
\int_0^t
\Vert \nabx(\overline{\bu}-\bv)\Vert_{L^2(\Omega_{\zeta})}^2\dt'
+ 
\gamma 
\int_0^t\Vert  \partial_{t'}\naby(\eta-\zeta) \Vert_{L^{2}(\omega )}^2
\dt'  
\\
\leq&K
\bigg[ 
 \Vert   \eta_\star-\zeta_\star \Vert_{L^{2}(\omega )}^2
+
\Vert  \naby^2(\eta_0-\zeta_0) \Vert_{L^{2}(\omega )}^2  
+
\Vert \overline{\bu}_0-\bv_0\Vert_{L^2(\Omega_{\zeta_0})}^2
\\&
\qquad+
(1+T)\Vert \overline{\rho}_0-\sigma_0\Vert_{L^2(\Omega_{\zeta_0})}^2  
+
\Vert \overline{\bT}_0-\bU_{0}\Vert_{L^2(\Omega_{\zeta_0})}^2
+
\varepsilon^2(1+T) 
\bigg].
\end{aligned}
\end{equation}
for all $t\in I$ where $(\overline{\bu}, \overline{p}, \overline{\rho},\overline{\bT})$ is given by \eqref{transHan} and where
\begin{align*}
K:=&C\,
e^{
c\int_0^t
 \left( 1+c(\gamma)
+
\Vert \partial_{t'}\overline{\bu}  \Vert_{L^2(\Omega_{\zeta})}^2
+
\Vert \overline{\bu}  \Vert_{W^{2,2}(\Ozeta)}^2
+\Vert \overline{p}  \Vert_{W^{1,2}(\Ozeta)}^2 
+
(1+\varepsilon^2)\Vert  \overline{\bT} \Vert_{W^{2,2}(\Omega_{\zeta})}^2
+
\Vert  \bU \Vert_{L^{\infty}(\Omega_{\zeta})}^2 
\right)\dt' }
\\&\times
e^{\left(
c\int_0^t
\left(  
(1+T)\Vert  \overline{\rho} \Vert_{W^{2,2}(\Omega_{\zeta})}^2
+
(1+T)\Vert \partial_{t'}\overline{\rho} \Vert_{L^{2}( \Omega_{\zeta})}^2 
+
\Vert \partial_{t'}\overline{\bT} \Vert_{L^{2}( \Omega_{\zeta})}^2   
\right)
\dt'\right)
}
\end{align*}
with $c$ and $C$ being constants independent of   $\varepsilon$.
\end{proposition}
\begin{proof}
As shown \cite[Section 4.1]{breit2022regularity} (see also \cite[Section 5]{BMSS}) for the solvent-structure subproblem and in \cite[Section 5.3]{mensah2023weak} for the solute subproblem, the transformed solution  $(\eta , \overline{\bu} , \overline{p} ,\overline{\rho} ,\overline{\bT} )$  solves
\begin{align}
\mathbb{B}_{\eta-\zeta}^\top : \nabx\overline{\bu}  = 0, \label{divfree2}
\\
\partial_t \overline{\rho} + (\overline{\bu} \cdot \nabx) \overline{\rho} 
= 
\varepsilon\Delx \overline{\rho} 
-
\varepsilon
\divx(\nabx \overline{\rho} (\mathbb{I}-\mathbb{A}_{\eta-\zeta}))
+
h_{\eta-\zeta}(\overline{\rho} ,\overline{\bu} )
,\label{rhoEqu2}
\\
\partial_t \overline{\bu}   = \Delx\overline{\bu}   
-
\nabx\overline{p} + \divx\overline{\bT} 
+
 \mathbf{h}_{\eta-\zeta}(\overline{\bu} )
-
\divx   \mathbb{G}_{\eta-\zeta}(\overline{\bu} ,\overline{p} ,\overline{\bT} ), \label{momEq2}
\\
\partial_t^2  \eta - \gamma\partial_t\partial_y^2  \eta +  \partial_y^4 \eta 
=  -([
\mathbb{A}_{\eta-\zeta}\nabx\overline{\bu}   
-
\mathbb{B}_{\eta-\zeta}(\overline{p} -\overline{\bT} )]\bn_{\zeta})\circ  \bm{\varphi}_{\zeta}\cdot\bn \,\det(\partial_y\bm{\varphi}_\zeta), 
\label{shellEQ2}
\\
\partial_t\overline{\bT} 
+
\overline{\bu}  \cdot\nabx \overline{\bT}  
=
\mathbb{W}(\nabx\overline{\bu}) \overline{\bT} 
+
\overline{\bT} \mathbb{W}( (\nabx\overline{\bu} )^\top)
-
2
(\overline{\bT} -\overline{\rho} \mathbb{I})
+
\varepsilon\Delx\overline{\bT} 
\nonumber
\\
-\varepsilon
\divx(\nabx \overline{\bT} (\mathbb{I}-\mathbb{A}_{\eta-\zeta}))
+
\mathbb{H}_{\eta-\zeta}(\overline{\rho} ,\overline{\bu} ,\overline{\bT} ) \label{solute2}
\end{align}
on $I \times \Omega_{\zeta}$ (with \eqref{shellEQ2} posed on $I\times\omega$) subject to the interface condition $\overline{\bu}\circ\bm{\varphi}_\zeta=(\partial_t\eta)\bn $ on $I\times \omega$ where
\begin{align*}
&J_{\eta-\zeta} :=\det(\nabx \bm{\Psi}_{\eta-\zeta}),
\\
&
\mathbb{B}_{\eta-\zeta}:=
J_{\eta-\zeta} \nabx \bm{\Psi}_{\eta-\zeta}^{-1}\circ \bm{\Psi}_{\eta-\zeta},
%\\
%&g_{\eta-\zeta}(\overline{\bu} ),
%=
%\nabx\overline{\bu} :(\mathbb{I}  -\mathbb{B}_{\eta-\zeta})
\\
&\mathbb{A}_{\eta-\zeta}:=
\mathbb{B}_{\eta-\zeta}(\nabx \bm{\Psi}_{\eta-\zeta}^{-1}\circ \bm{\Psi}_{\eta-\zeta})^\top,
\\&
\mathbb{H}_{\eta-\zeta}(\overline{\rho} ,\overline{\bu} ,\overline{\bT} )
=
(1-J_{\eta-\zeta})\partial_t\overline{\bT} 
-
J_{\eta-\zeta} \nabx \overline{\bT} \cdot\partial_t\bm{\Psi}_{\eta-\zeta}^{-1}\circ \bm{\Psi}_{\eta-\zeta}
+
\mathbb{W}(\nabx\overline{\bu}) (\mathbb{B}_{\eta-\zeta}-\mathbb{I})\overline{\bT} 
\\&
\qquad\qquad
+
\overline{\bT} (\mathbb{B}_{\eta-\zeta}-\mathbb{I})^\top\mathbb{W}( (\nabx\overline{\bu} )^\top)
+
\overline{\bu}  \cdot\nabx \overline{\bT}  (\mathbb{I}-\mathbb{B}_{\eta-\zeta})
+
2
(1-J_{\eta-\zeta}) (\overline{\bT} -\overline{\rho} \mathbb{I}),
\\&
\mathbb{G}_{\eta-\zeta}(\overline{\bu} ,\overline{p} ,\overline{\bT} )
=
(\mathbb{I} -\mathbb{A}_{\eta-\zeta}) \nabx\overline{\bu} 
-
(\mathbb{I} -\mathbb{B}_{\eta-\zeta})(\overline{p} -\overline{\bT} ),
\\&
 \mathbf{h}_{\eta-\zeta}(\overline{\bu} )
 =
 (1-J_{\eta-\zeta})\partial_t\overline{\bu} 
-
J_{\eta-\zeta} \nabx \overline{\bu} \partial_t\bm{\Psi}_{\eta-\zeta}^{-1}\circ \bm{\Psi}_{\eta-\zeta}
-
\nabx\overline{\bu} \mathbb{B}_{\eta-\zeta}\overline{\bu} ,
\\&
h_{\eta-\zeta}(\overline{\rho} ,\overline{\bu} )
=
(1-J_{\eta-\zeta})\partial_t\overline{\rho} 
-
J_{\eta-\zeta} \nabx \overline{\rho} \cdot\partial_t\bm{\Psi}_{\eta-\zeta}^{-1}\circ \bm{\Psi}_{\eta-\zeta}
+
\overline{\bu}  \cdot\nabx \overline{\rho}  (\mathbb{I}-\mathbb{B}_{\eta-\zeta}).
\end{align*} 
Furthermore, by following a similar argument as \cite[Lemma 4.2]{breit2022regularity}, we have that $(\eta, \overline{\bu}, \overline{p},\overline{\rho},\overline{\bT} )$ is a strong solution in the sense that it inherits on $\Ozeta$, the same regularity properties of $(\eta, \bu, p, \rho,\bT)$ in Definition \ref{def:strongSolution}. With this preparation in hand,
we are now going to obtain estimates for the differences
$
(\eta-\zeta ,  \overline{\bu}-\bv, \overline{p}-\pi,  \overline{\rho}-\sigma, \overline{\bT} -\bU)
$ where we start with the components that describe the solvent-structure subproblem.
\subsection{Estimate for the solvent-structure interaction}
For the structure, we aim to obtain an estimate for the following sum
\begin{align*}
 \Vert  \partial_t(\eta-\zeta)(t) \Vert_{L^{2}(\omega )}^2
+
\Vert  \naby^2(\eta-\zeta)(t) \Vert_{L^{2}(\omega )}^2  
+
\gamma
\int_0^t
\Vert  \partial_t\naby(\eta-\zeta) \Vert_{L^{2}(\omega )}^2\dt'
\end{align*}
for any $t\in I$.  In this regard, given the basic identity $\tfrac{1}{2}\vert a -b\vert^2= \tfrac{1}{2}\vert a \vert^2+\tfrac{1}{2}\vert  b\vert^2-ab$ which holds for real numbers and the fact that the equations (the first two equations below follows by testing \eqref{shellEQ2} and \eqref{shellEQL} with $\partial_t\eta$ and $\partial_t\zeta$, respectively)
\begin{align*}
\frac{1}{2}\frac{\dd}{\dt}\Vert \partial_t\eta(t)\Vert_{L^2(\omega)}^2
=&-
\frac{1}{2}\frac{\dd}{\dt}\Vert \partial_y^2\eta(t)\Vert_{L^2(\omega)}^2
-\gamma\Vert\partial_t\partial_y\eta\Vert_{L^2(\omega)}^2
\\&-\int_\omega\partial_t\eta([
\mathbb{A}_{\eta-\zeta}\nabx\overline{\bu}   
-
\mathbb{B}_{\eta-\zeta}(\overline{p} -\overline{\bT} )]\bn_{\zeta})\circ  \bm{\varphi}_{\zeta}\cdot\bn \,\det(\partial_y\bm{\varphi}_\zeta)\dy
\\
\frac{1}{2}\frac{\dd}{\dt}\Vert \partial_t\zeta(t)\Vert_{L^2(\omega)}^2
=&
-
\frac{1}{2}\frac{\dd}{\dt}\Vert \partial_y^2\zeta(t)\Vert_{L^2(\omega)}^2
-\gamma\Vert\partial_t\partial_y\zeta\Vert_{L^2(\omega)}^2
\\&
-\int_\omega\partial_t\zeta\,( \mathbb{S}_*\bn_\zeta )\circ \bm{\varphi}_\zeta\cdot\bn \,\det(\partial_y\bm{\varphi}_\zeta) \dy
\\
-\frac{\dd}{\dt}
\int_{\omega}\partial_t\zeta\partial_t\eta\dy
=&
-\int_{\omega}\partial_t^2\zeta\partial_t\eta\dy
-\int_{\omega}\partial_t\zeta\partial_t^2\eta\dy
%\\&= 
%\int_\omega
%\partial_t\zeta
% \big[\partial_y^4 \eta 
% - \gamma\partial_t\partial_y^2  \eta
%+([
%\mathbb{A}_{\eta-\zeta}\nabx\overline{\bu}   
%-
%\mathbb{B}_{\eta-\zeta}(\overline{p} -\overline{\bT} )]\bn_{\zeta})\circ  \bm{\varphi}_{\zeta}\cdot\bn \,\det(\partial_y\bm{\varphi}_\zeta)
%\big]\dy
%\\&+ 
%\int_\omega
%\partial_t\eta
% \big[\partial_y^4 \zeta  
%+( \mathbb{S}_*\bn_\zeta )\circ \bm{\varphi}_\zeta\cdot\bn \,\det(\partial_y\bm{\varphi}_\zeta)
%\big]\dy
%\\&= 
%\int_\omega
% \big[\partial_t\partial_y^2\zeta\partial_y^2 \eta 
% - \gamma\partial_t\partial_y\zeta\partial_t\partial_y \eta]
%+\partial_t\zeta
% \big[([
%\mathbb{A}_{\eta-\zeta}\nabx\overline{\bu}   
%-
%\mathbb{B}_{\eta-\zeta}(\overline{p} -\overline{\bT} )]\bn_{\zeta})\circ  \bm{\varphi}_{\zeta}\cdot\bn \,\det(\partial_y\bm{\varphi}_\zeta)
%\big]\dy
%\\&+ 
%\int_\omega
% \big[\partial_t\partial_y^2\eta\partial_y^2 \zeta  
%+\partial_t\eta( \mathbb{S}_*\bn_\zeta )\circ \bm{\varphi}_\zeta\cdot\bn \,\det(\partial_y\bm{\varphi}_\zeta)
%\big]\dy
\\=& 
\frac{\dd}{\dt}
\int_\omega
\partial_y^2\zeta\partial_y^2 \eta \dy
+
2\gamma
\int_\omega 
 \partial_t\partial_y\zeta\partial_t\partial_y \eta
 \dy
\\&+ 
\int_\omega
\partial_t\eta( \mathbb{S}_*\bn_\zeta )\circ \bm{\varphi}_\zeta\cdot\bn \,\det(\partial_y\bm{\varphi}_\zeta)
\dy
\\&+
 \int_\omega
 \partial_t\zeta
([
\mathbb{A}_{\eta-\zeta}\nabx\overline{\bu}   
-
\mathbb{B}_{\eta-\zeta}(\overline{p} -\overline{\bT} )]\bn_{\zeta})\circ  \bm{\varphi}_{\zeta}\cdot\bn \,\det(\partial_y\bm{\varphi}_\zeta)
\dy
\end{align*}
follow from \eqref{shellEQ2} and \eqref{shellEQL}, we obtain
%\begin{align*}
%\partial_t^2 (\eta-\zeta)  +  \partial_y^4 (\eta-\zeta) 
%= \gamma\partial_t\partial_y^2  \eta
%+
%( \mathbb{S}_*\bn_\zeta )\circ \bm{\varphi}_\zeta\cdot\bn \,\det(\partial_y\bm{\varphi}_\zeta)
% -
% ([
%\mathbb{A}_{\eta-\zeta}\nabx\overline{\bu}   
%-
%\mathbb{B}_{\eta-\zeta}(\overline{p} -\overline{\bT} )]\bn_{\zeta})\circ  \bm{\varphi}_{\zeta})\cdot\bn\,\det(\partial_y\bm{\varphi}_\zeta)
%\end{align*}
%on $I\times\omega$. Testing the above with $\partial_t (\eta-\zeta)$ yields
\begin{equation}
\begin{aligned}
\label{shellDiv1}
\frac{1}{2}\frac{\dd}{\dt}\big( \Vert&  \partial_t(\eta-\zeta)(t) \Vert_{L^{2}(\omega )}^2
+
\Vert  \naby^2(\eta-\zeta)(t) \Vert_{L^{2}(\omega )}^2 \big)
+
\gamma\Vert  \partial_t\naby(\eta-\zeta) \Vert_{L^{2}(\omega )}^2
\\=&
\int_{\omega}( \mathbb{S}_*\bn_\zeta )\circ \bm{\varphi}_\zeta\cdot(\partial_t(\eta-\zeta)\bn) \det(\partial_y\bm{\varphi}_\zeta)
\dy
\\&-
\int_{\omega}([
\mathbb{A}_{\eta-\zeta}\nabx\overline{\bu}   
-
\mathbb{B}_{\eta-\zeta}(\overline{p} -\overline{\bT} )]\bn_{\zeta})\circ  \bm{\varphi}_{\zeta})\cdot(\partial_t(\eta-\zeta)\bn)  
\det(\partial_y\bm{\varphi}_\zeta)
\dy.
\end{aligned}
\end{equation}
In our next goal, we obtain a difference estimate for the equation of the solvent. More precisely, we seek to obtain an estimate for
\begin{align*}
 \Vert (\overline{\bu}-\bv)(t)\Vert_{L^2(\Omega_{\zeta})}^2
+
\int_0^t
\Vert \nabx(\overline{\bu}-\bv)\Vert_{L^2(\Omega_{\zeta})}^2\dt'
\end{align*} 
for any $t\in I$.
Here, because $\overline{\bu}$ is not solenoidal but $\bv$ is, we need to be extra careful. In particular, $\overline{\bu}$ cannot be a test function for the evolution equation of $\bv$ since the former does not have the solenoidal property that the latter possesses. As we shall soon see, the key to remedying this problem is a suitable addition and subtraction of the Bogovskij operator introduced in Section \ref{sec:prelim} to render test functions solenoidal. Related to this, the introduction of the Bogovskij operator kills the undesirable pressure $\pi$ in the evolution of $\bv$ once $\bv$ interacts with the evolution of $\overline{\bu}$ via the difference $\Vert\overline{\bu}-\bv\Vert_{L^2(\Ozeta)}^2$.
Again, we rely on the basic identity $\tfrac{1}{2}\vert a -b\vert^2= \tfrac{1}{2}\vert a \vert^2+\tfrac{1}{2}\vert  b\vert^2-ab$ to obtain this difference estimate. In this regard, we first observe that the energy identity
\begin{align*}
\frac{1}{2}\frac{\dd}{\dt} \Vert  \bv(t) \Vert_{L^{2}( \Omega_{\zeta})}^2
+
\Vert\nabx \bv\Vert_{L^{2}( \Omega_{\zeta})}^2
=&\int_\omega ( \mathbb{S}_*\bn_\zeta )\circ \bm{\varphi}_\zeta\cdot(\partial_t\zeta\bn) \,\det(\partial_y\bm{\varphi}_\zeta)\dy
\\&-
\int_{\Ozeta}\bU:\nabx\bv\dx
\end{align*}
holds by testing \eqref{momEqL} with $\bv$ and using the interface condition \eqref{interfaceL}. 
Next, by applying Reynold's transport theorem to the integral $\int_{\Ozeta}\overline{\bu}\cdot\partial_t\overline{\bu}\dx$ and
using the interface condition, we obtain the identity
\begin{align*}
\frac{1}{2}\frac{\dd}{\dt} \Vert  \overline{\bu}(t) \Vert_{L^{2}( \Omega_{\zeta})}^2 
+
\Vert\nabx \overline{\bu}\Vert_{L^{2}( \Omega_{\zeta})}^2
=&
\frac{1}{2}
\int_{\partial \Omega_{\zeta}}(\bn\partial_t\zeta)\circ\bm{\varphi}_\zeta^{-1}\cdot\bn_\zeta\vert\overline{\bu}\vert^2\dd\mathcal{H}^1
\\&+
\int_{\Ozeta}\overline{\bu}\cdot\partial_t\overline{\bu}\dx
+
\int_{\Ozeta}\nabx \overline{\bu}:\nabx\overline{\bu}\dx
\end{align*}
without any recourse to \eqref{momEq2}.
Similarly, the Reynold's transport theorem applied to the product yield
\begin{align*}
- \frac{\dd}{\dt}\int_{\Omega_{\zeta}}\overline{\bu}\cdot\bv\dx 
=&
-\int_{\partial \Omega_{\zeta}}(\bn\partial_t\zeta)\circ\bm{\varphi}_\zeta^{-1}\cdot\bn_\zeta\overline{\bu}\cdot\bv\dd\mathcal{H}^1
\\&-
\int_{\Omega_{\zeta}}(\overline{\bu}-\mathrm{Bog}_{\zeta}(\divx \overline{\bu}))\cdot\partial_t\bv\dx
-
\int_{\Ozeta}\overline{\bu}\cdot\partial_t\overline{\bu}\dx
\\&-
\int_{\Omega_{\zeta}} (\bv -\overline{\bu} + \mathrm{Bog}_{\zeta}(\divx \overline{\bu}) )\cdot\partial_t\overline{\bu}\dx
\\&+
\int_{\Omega_{\zeta}} \mathrm{Bog}_{\zeta}(\divx \overline{\bu})\cdot\partial_t(\overline{\bu} - \bv)\dx.
\end{align*}
We note that the addition and subtraction of the Bogovskij operator (i.e. an addition of zero) in the last equation above ensures that any time derivative of a function is being tested with a function with the same or better properties.
Moving on, by recalling the basic identity $\tfrac{1}{2}\vert a -b\vert^2= \tfrac{1}{2}\vert a \vert^2+\tfrac{1}{2}\vert  b\vert^2-ab$, it follows that
\begin{equation}
\begin{aligned}
\label{velDiff1}
\frac{1}{2}\frac{\dd}{\dt}\Vert (\overline{\bu}-\bv)(t)\Vert_{L^2(\Omega_{\zeta})}^2
&+ 
\Vert \nabx(\overline{\bu}-\bv)\Vert_{L^2(\Omega_{\zeta})}^2
\\=& 
\int_\omega ( \mathbb{S}_*\bn_\zeta )\circ \bm{\varphi}_\zeta\cdot (\partial_t \zeta \bn) \,\det(\partial_y\bm{\varphi}_\zeta)\dy 
-
\int_{\Ozeta}\bU:\nabx\bv\dx
\\
&+
\frac{1}{2}
\int_{\partial \Omega_{\zeta}}(\bn\partial_t\zeta)\circ\bm{\varphi}_\zeta^{-1}\cdot\bn_\zeta(\vert\overline{\bu}\vert^2 - 2 \overline{\bu}\cdot\bv)\dd\mathcal{H}^1 
\\&+
\int_{\Ozeta}\nabx \overline{\bu}:\nabx\overline{\bu}\dx
-2
\int_{\Ozeta}\nabx \bv:\nabx\overline{\bu}\dx
\\&
-
\int_{\Omega_{\zeta}}(\overline{\bu}-\mathrm{Bog}_{\zeta}(\divx \overline{\bu}))\cdot\partial_t\bv\dx 
\\&-
\int_{\Omega_{\zeta}} (\bv -\overline{\bu} + \mathrm{Bog}_{\zeta}(\divx \overline{\bu}) \cdot\partial_t\overline{\bu}\dx
\\&+
\int_{\Omega_{\zeta}} \mathrm{Bog}_{\zeta}(\divx \overline{\bu})\cdot\partial_t(\overline{\bu} - \bv)\dx.
\end{aligned}
\end{equation}
We are now going to rewrite any time-derivative of a velocity term on the right-hand side of the identity above by either substituting in \eqref{momEqL}, \eqref{momEq2}, or performing an integration by part argument.  
We do this because, in particular, $\partial_t\bv$ does not possess the required \textit{explicit} regularity according to Definition \ref{def:weaksolmartFPBounded}. As would be seen after substituting in \eqref{momEqL}, this deficiency in regularity is only explicit. Implicitly, however, $\partial_t\bv$ belong to some sufficient negative Sobolev space because all terms to the right of the equation $\partial_t\bv=\ldots$ (i.e. \eqref{momEqL}) also belong to at least a negative Sobolev space according to Definition \ref{def:weaksolmartFPBounded}. 

Indeed, by using \eqref{momEqL}, integrating by parts where necessary, and using the properties of the Bogovskij operator, we note that
\begin{align*}
-
\int_{\Omega_{\zeta}}(\overline{\bu}-\mathrm{Bog}_{\zeta}(\divx \overline{\bu}))\cdot\partial_t\bv\dx
=&
%\int_{\Omega_{\zeta}}(\overline{\bu}-\mathrm{Bog}_{\zeta}(\divx \overline{\bu}))\cdot[ (\bv\cdot \nabx)\bv-\delx \bv +\nabx \pi 
%-
%\divx   \bU]\dx
%\\&
%=
\int_{\Omega_{\zeta}}(\overline{\bu}-\mathrm{Bog}_{\zeta}(\divx \overline{\bu}))\cdot ((\bv\cdot \nabx)\bv)\dx
\\&
-
\int_\omega (  \mathbb{S}_* \bn_\zeta )\circ \bm{\varphi}_\zeta\cdot(\partial_t\eta\bn) \,\det(\partial_y\bm{\varphi}_\zeta)\dy
\\&+
\int_{\Omega_{\zeta}}\nabx(\overline{\bu} -\mathrm{Bog}_{\zeta}(\divx \overline{\bu})): \nabx \bv\dx
\\&+
\int_{\Omega_{\zeta}}\nabx(\overline{\bu}-\mathrm{Bog}_{\zeta}(\divx \overline{\bu})):   \bU \dx
\end{align*} 
and since the Bogovskij operator has zero trace, we also have
\begin{align*}
\int_{\Omega_{\zeta}} \mathrm{Bog}_{\zeta}(\divx \overline{\bu})\cdot\partial_t(\overline{\bu} - \bv)\dx 
=&
\frac{\dd}{\dt}
\int_{\Omega_{\zeta}} \mathrm{Bog}_{\zeta}(\divx \overline{\bu})\cdot(\overline{\bu} - \bv)\dx 
\\&-
\int_{\Omega_{\zeta}} \partial_t\mathrm{Bog}_{\zeta}(\divx \overline{\bu})\cdot(\overline{\bu} - \bv)\dx.
\end{align*}
Now note that the identity 
\begin{align*}
\int_{\Omega_{\zeta}}\nabx(\bv -\overline{\bu} + \mathrm{Bog}_{\zeta}(\divx \overline{\bu})): \overline{p}\mathbb{I} 
\dx
=
\int_{\Omega_{\zeta}}\divx(\bv -\overline{\bu} + \mathrm{Bog}_{\zeta}(\divx \overline{\bu}))\overline{p}
\dx=0
\end{align*}
hold because of \eqref{divfreeL} and the fact that the operator $\mathrm{Bog}_{\zeta}$ serves as the inverse of the divergence operator. Also, we can use the identities
$(\nabx\overline{\bu})\overline{\bu}=(\nabx\overline{\bu})\mathbb{I}\overline{\bu}=(\overline{\bu}\cdot\nabx)\overline{\bu}$
to rewrite $ \mathbf{h}_{\eta-\zeta}(\overline{\bu} )$ in \eqref{momEq2} as
\begin{align*}
 \mathbf{h}_{\eta-\zeta}(\overline{\bu} )
 =
 (1-J_{\eta-\zeta})\partial_t\overline{\bu} 
-
J_{\eta-\zeta} \nabx \overline{\bu} \partial_t\bm{\Psi}_{\eta-\zeta}^{-1}\circ \bm{\Psi}_{\eta-\zeta}
+
\nabx\overline{\bu}(\mathbb{I}- \mathbb{B}_{\eta-\zeta})\overline{\bu}
-(\overline{\bu}\cdot\nabx)\overline{\bu}.
\end{align*}
Consequently, by using the trace-free property of the Bogovskij operator, we obtain
\begin{align*}
-
\int_{\Omega_{\zeta}} &(\bv -\overline{\bu} + \mathrm{Bog}_{\zeta}(\divx \overline{\bu}) \cdot\partial_t\overline{\bu}\dx
 =  
\int_{\Omega_{\zeta}}\nabx (\bv -\overline{\bu} + \mathrm{Bog}_{\zeta}(\divx \overline{\bu}))  :\nabx\overline{\bu}   
\dx
\\&+
\int_{\Omega_{\zeta}}\nabx (\bv -\overline{\bu} + \mathrm{Bog}_{\zeta}(\divx \overline{\bu}))  :\overline{\bT}  
\dx
\\&+
\int_{\Omega_{\zeta}}  (\bv -\overline{\bu} + \mathrm{Bog}_{\zeta}(\divx \overline{\bu}))  \cdot 
((\overline{\bu}\cdot\nabx)\overline{\bu})
\dx
\\&-
\int_{\Omega_{\zeta}}  (\bv -\overline{\bu} + \mathrm{Bog}_{\zeta}(\divx \overline{\bu}))  \cdot 
[ (1-J_{\eta-\zeta})\partial_t\overline{\bu} 
-
J_{\eta-\zeta} \nabx \overline{\bu} \partial_t\bm{\Psi}_{\eta-\zeta}^{-1}\circ \bm{\Psi}_{\eta-\zeta}]
\dx
\\&-
\int_{\Omega_{\zeta}}  (\bv -\overline{\bu} + \mathrm{Bog}_{\zeta}(\divx \overline{\bu}))  \cdot 
[
\nabx\overline{\bu}(\mathbb{I}- \mathbb{B}_{\eta-\zeta})\overline{\bu}]
\dx
\\&-
\int_{\Omega_{\zeta}}\nabx (\bv -\overline{\bu} + \mathrm{Bog}_{\zeta}(\divx \overline{\bu}))  :[(\mathbb{I} -\mathbb{A}_{\eta-\zeta}) \nabx\overline{\bu} 
-
(\mathbb{I} -\mathbb{B}_{\eta-\zeta})(\overline{p} -\overline{\bT} )]
\dx
\\&+
\int_{\omega}([
\mathbb{A}_{\eta-\zeta}\nabx\overline{\bu}   
-
\mathbb{B}_{\eta-\zeta}(\overline{p} -\overline{\bT} )]\bn_{\zeta})\circ  \bm{\varphi}_{\zeta})\cdot(\partial_t(\eta-\zeta)\bn)  
\det(\partial_y\bm{\varphi}_\zeta)
\dy.
\end{align*}
If we substitute the identities above into \eqref{velDiff1} and  sum the resulting equation with \eqref{shellDiv1},
then we obtain the equation
\begin{equation}
\begin{aligned}
\label{velShellDiff1}
\frac{1}{2}\frac{\dd}{\dt}\big(\Vert (\overline{\bu}&-\bv)(t)\Vert_{L^2(\Omega_{\zeta})}^2
+
 \Vert  \partial_t(\eta-\zeta)(t) \Vert_{L^{2}(\omega )}^2
+
\Vert  \naby^2(\eta-\zeta)(t) \Vert_{L^{2}(\omega )}^2 \big)
\\&\qquad+ 
\Vert \nabx(\overline{\bu}-\bv)\Vert_{L^2(\Omega_{\zeta})}^2
+
\gamma \Vert  \partial_t\naby(\eta-\zeta) \Vert_{L^{2}(\omega )}^2
\\
=& 
\frac{1}{2}
\int_{\partial \Omega_{\zeta}}(\bn\partial_t\zeta)\circ\bm{\varphi}_\zeta^{-1}\cdot\bn_\zeta(\vert\overline{\bu}\vert^2 - 2 \overline{\bu}\cdot\bv)\dd\mathcal{H}^1 
\\&+
\int_{\Omega_{\zeta}}(\overline{\bu}-\mathrm{Bog}_{\zeta}(\divx \overline{\bu}))\cdot ((\bv\cdot \nabx)\bv)\dx 
\\&+
\int_{\Omega_{\zeta}}  (\bv -\overline{\bu} + \mathrm{Bog}_{\zeta}(\divx \overline{\bu}))  \cdot 
((\overline{\bu}\cdot\nabx)\overline{\bu})
\dx
\\&+
\int_{\Omega_{\zeta}}\nabx\mathrm{Bog}_{\zeta}(\divx \overline{\bu}): \nabx (\overline{\bu} - \bv)\dx 
\\&+
\int_{\Omega_{\zeta}}\nabx (\bv -\overline{\bu} + \mathrm{Bog}_{\zeta}(\divx \overline{\bu}))  :(\overline{\bT}  -\bU)
\dx
\\&-
\int_{\Omega_{\zeta}}  (\bv -\overline{\bu} + \mathrm{Bog}_{\zeta}(\divx \overline{\bu}))  \cdot 
(1-J_{\eta-\zeta})\partial_t\overline{\bu} 
\dx
\\&-
\int_{\Omega_{\zeta}}  (\bv -\overline{\bu} + \mathrm{Bog}_{\zeta}(\divx \overline{\bu}))  \cdot 
\nabx\overline{\bu}(\mathbb{I}- \mathbb{B}_{\eta-\zeta})\overline{\bu}
\dx
\\&+
\int_{\Omega_{\zeta}}  (\bv -\overline{\bu} + \mathrm{Bog}_{\zeta}(\divx \overline{\bu}))  \cdot 
J_{\eta-\zeta} \nabx \overline{\bu} \partial_t\bm{\Psi}_{\eta-\zeta}^{-1}\circ \bm{\Psi}_{\eta-\zeta}
\dx
\\&-
\int_{\Omega_{\zeta}}\nabx (\bv -\overline{\bu} + \mathrm{Bog}_{\zeta}(\divx \overline{\bu}))  :[(\mathbb{I} -\mathbb{A}_{\eta-\zeta}) \nabx\overline{\bu} ]
\dx 
\\&+
\int_{\Omega_{\zeta}}\nabx (\bv -\overline{\bu} + \mathrm{Bog}_{\zeta}(\divx \overline{\bu}))  :[
(\mathbb{I} -\mathbb{B}_{\eta-\zeta})(\overline{p} -\overline{\bT} )]
\dx 
\\&+
\frac{\dd}{\dt}
\int_{\Omega_{\zeta}} \mathrm{Bog}_{\zeta}(\divx \overline{\bu})\cdot(\overline{\bu} - \bv)\dx 
-
\int_{\Omega_{\zeta}} \partial_t\mathrm{Bog}_{\zeta}(\divx \overline{\bu})\cdot(\overline{\bu} - \bv)\dx
\\&
=: J_1+J_2+\ldots+J_{12}.
\end{aligned}
\end{equation}
We are now going to estimate the terms on the right. First of all, we observe that since $\overline{\bu}-\mathrm{Bog}_{\zeta}(\divx \overline{\bu})$ is divergence-free and the Bogovskij operator has zero trace, we can use the interface condition  \eqref{interfaceL} to obtain
\begin{align*}
J_2
&=
\frac{1}{2}
\int_{\Omega_{\zeta}}(\overline{\bu}-\mathrm{Bog}_{\zeta}(\divx \overline{\bu}))\cdot  \nabx\vert\bv\vert^2\dx
\\&
=
\frac{1}{2}
\int_{\partial\Omega_{\zeta}}(\bn\partial_t\eta)\circ\bm{\varphi}_\zeta^{-1}\cdot  \bn_\zeta\vert\bv\vert^2\dd\mathcal{H}^1
\\
&=
\frac{1}{2}
\int_{\omega}(\bn\partial_t\eta) \cdot\bn_\zeta\circ\bm{\varphi}_\zeta((\bn\partial_t\zeta) \cdot(\bn\partial_t\zeta) )\,\det(\partial_y\bm{\varphi}_\zeta)\dy 
\\
&=
\frac{1}{2}
\int_{\omega}(\bn\partial_t\zeta) \cdot\bn_\zeta\circ\bm{\varphi}_\zeta((\bn\partial_t\eta) \cdot(\bn\partial_t\zeta) )\,\det(\partial_y\bm{\varphi}_\zeta)\dy 
\\
&=
\frac{1}{2}
\int_{\partial \Omega_{\zeta}}(\bn\partial_t\zeta)\circ\bm{\varphi}_\zeta^{-1}\cdot\bn_\zeta\overline{\bu}\cdot\bv\dd\mathcal{H}^1
\end{align*}
and similarly,
\begin{align*}
J_3
&=
\frac{1}{2}
\int_{\partial\Omega_{\zeta}}(\bn\partial_t(\zeta - \eta))\circ\bm{\varphi}_\zeta^{-1}\cdot  \bn_\zeta\vert\overline{\bu}\vert^2\dd\mathcal{H}^1
\\&
=
\frac{1}{2}
\int_{\partial\Omega_{\zeta}}(\bn\partial_t\eta)\circ\bm{\varphi}_\zeta^{-1}\cdot  \bn_\zeta\overline{\bu}\cdot\bv\dd\mathcal{H}^1
-
\frac{1}{2}
\int_{\partial\Omega_{\zeta}}(\bn\partial_t\eta)\circ\bm{\varphi}_\zeta^{-1}\cdot  \bn_\zeta\vert\overline{\bu}\vert^2\dd\mathcal{H}^1.
\end{align*}
Thus, it follows from the Trace Theorem that
\begin{align*}
J_1+J_2+J_3
&=-
\frac{1}{2}
\int_{\partial\Omega_{\zeta}}(\bn\partial_t(\zeta - \eta))\circ\bm{\varphi}_\zeta^{-1}\cdot  \bn_\zeta\overline{\bu}\cdot(\overline{\bu}-\bv)\dd\mathcal{H}^1
\\&\leq
\delta
\Vert \nabx (\overline{\bu} - \bv)\Vert_{L^2(\Ozeta)}^2
+
c(\delta)
  \Vert \partial_t(\eta-\zeta)\Vert_{L^{2}(\omega)}^2
   \Vert \overline{\bu} \Vert_{W^{2,2}(\Ozeta)}^2
\end{align*}
holds for any $\delta>0$.
Now since the identity $\divx \overline{\bu}=\mathbb{I}:\nabx \overline{\bu}=(\mathbb{I}-\mathbb{B}_{\eta-\zeta}^\top):\nabx \overline{\bu}$ follows from \eqref{divfree2} and the continuity property of the Bogovskij operator yields
\begin{align*}
\Vert \nabx\mathrm{Bog}_{\zeta}(\divx \overline{\bu})\Vert_{L^2(\Ozeta)}
 \lesssim 
 \Vert \divx \overline{\bu} \Vert_{L^2(\Ozeta)},
\end{align*}
it follows that
\begin{align*} 
J_4&
\lesssim
 \Vert \nabx (\overline{\bu} - \bv)\Vert_{L^2(\Ozeta)}
  \Vert \naby (\eta-\zeta)\Vert_{L^4(\omega)}
   \Vert \nabx \overline{\bu} \Vert_{L^4(\Ozeta)}
\\&\leq
\delta
\Vert \nabx (\overline{\bu} - \bv)\Vert_{L^2(\Ozeta)}^2
+
c(\delta)
  \Vert \eta-\zeta\Vert_{W^{2,2}(\omega)}^2
   \Vert \overline{\bu} \Vert_{W^{2,2}(\Ozeta)}^2
\end{align*}
holds for any $\delta>0$ and similarly
\begin{align*} 
J_5&\leq
\frac{1}{2} 
\Vert  \overline{\bT} - \bU\Vert_{L^2(\Ozeta)}^2
+
\frac{1}{2}\Vert \nabx (\overline{\bu} - \bv)\Vert_{L^2(\Ozeta)}^2
\\&
+
\delta_*
\Vert  \overline{\bT} - \bU\Vert_{L^2(\Ozeta)}^2
+
c(\delta_*)
  \Vert \eta-\zeta\Vert_{W^{2,2}(\omega)}^2
   \Vert \overline{\bu} \Vert_{W^{2,2}(\Ozeta)}^2
\end{align*}
holds for any $\delta_*>0$. The subscript  $*$ is just a marker because we do not want to absorb this particular term in the left-hand side just yet. This will become clear below.
\\
Next,  by using the estimate
\begin{align*}
\Vert \mathrm{Bog}_{\zeta}(\divx \overline{\bu})\Vert_{W^{1,2}(\Ozeta)}
=
\Vert \mathrm{Bog}_{\zeta}(\divx (\overline{\bu}-\bv))\Vert_{W^{1,2}(\Ozeta)}
\lesssim
\Vert
\overline{\bu}-\bv \Vert_{W^{1,2}(\Ozeta)}
\end{align*}
we have
\begin{align*}
J_6&\leq
 \delta
 \Vert \nabx(\overline{\bu}-\bv) \Vert_{L^2(\Omega_{\zeta})}^2
 +
 c(\delta)\Vert \partial_t\overline{\bu}  \Vert_{L^2(\Omega_{\zeta})}^2
 \Vert \eta - \zeta \Vert_{W^{2,2}(\omega)}^2,
 \\
 J_7&\leq
 \delta
 \Vert \nabx(\overline{\bu}-\bv) \Vert_{L^2(\Omega_{\zeta})}^2
 +
 c(\delta)\Vert  \overline{\bu}  \Vert_{W^{2,2}(\Omega_{\zeta})}^2
 \Vert \eta - \zeta \Vert_{W^{2,2}(\omega)}^2
\end{align*} 
for any $\delta>0$ whereas by Ladyzhenskaya's inequality,
\begin{align*}
J_8&\lesssim
 \Vert  \overline{\bu}-\bv\Vert_{L^4(\Omega_{\zeta})}
 \Vert \nabx \overline{\bu}  \Vert_{L^{4}(\Omega_{\zeta})}
 \Vert \partial_t(\eta - \zeta) \Vert_{L^{2}(\omega)}
\\
&\leq
\delta
\Vert \nabx(\overline{\bu}-\bv) \Vert_{L^2(\Ozeta)}^2
+
c(\delta)
\Vert \overline{\bu} - \bv\Vert_{L^2(\Ozeta)}^2
+
c(\delta)
  \Vert \partial_t(\eta-\zeta)\Vert_{L^{2}(\omega)}^2
   \Vert \overline{\bu} \Vert_{W^{2,2}(\Ozeta)}^2.
\end{align*}
Similar to $J_6$, we also have
\begin{align*}
J_{9} 
\leq
\delta
\Vert \nabx(\overline{\bu}-\bv) \Vert_{L^2(\Ozeta)}^2
+
c(\delta)
\Vert \overline{\bu}  \Vert_{W^{2,2}(\Ozeta)}^2
\Vert \eta-\zeta\Vert_{W^{2,2}(\omega)}^2
\end{align*}
and
\begin{align*}
J_{10}
\leq
\delta
\Vert \nabx(\overline{\bu}-\bv) \Vert_{L^2(\Ozeta)}^2
+
c(\delta)
\big(\Vert \overline{p}  \Vert_{W^{1,2}(\Ozeta)}^2
+\Vert \overline{\bT}  \Vert_{W^{1,2}(\Ozeta)}^2
\big)\Vert \eta-\zeta\Vert_{W^{2,2}(\omega)}^2.
\end{align*}
Next, we observe that the identity $\divx \overline{\bu} =(\mathbb{I}-\mathbb{B}_{\eta-\zeta}^\top):\nabx \overline{\bu}$ and the continuous embedding  $L^{6/5}(\Ozeta)\hookrightarrow W^{-1,2}(\Ozeta)$ 
leads to
\begin{align*}
\Vert  \mathrm{Bog}_{\zeta}(\divx \overline{\bu})
\Vert_{L^2(\Ozeta)}^2
&\lesssim
\Vert  (\mathbb{I}-\mathbb{B}_{\eta-\zeta}^\top):\nabx \overline{\bu}
\Vert_{L^{6/5}(\Ozeta)}^2
\\
&\lesssim
\Vert  \naby(\eta-\zeta)\Vert_{L^3(\omega)}^2\Vert \nabx \overline{\bu}
\Vert_{L^2(\Ozeta)}^2
\\
&\lesssim
\Vert  \naby(\eta-\zeta)\Vert_{L^2(\omega)}^{4/3}\Vert  \naby^2(\eta-\zeta)\Vert_{L^2(\omega)}^{2/3}\Vert \nabx \overline{\bu}
\Vert_{L^2(\Ozeta)}^2
\\
&\leq\delta
\Vert  \naby^2(\eta-\zeta)\Vert_{L^2(\omega)}^{2}
+
c(\delta)\Vert  \naby(\eta-\zeta)\Vert_{L^2(\omega)}^{2}
\Vert \nabx \overline{\bu}
\Vert_{L^2(\Ozeta)}^3
\end{align*}
where for any $t\in I$, the $1$-dimensional Agmon's inequality in time yields
\begin{align*}
\Vert  \naby(\eta-\zeta)(t)\Vert_{L^2(\omega)}^{2}
&\leq
\sup_{t'\in[0,t)}\Vert  \naby(\eta-\zeta)(t')\Vert_{L^2(\omega)}^{2} 
\\
&\lesssim
\bigg(
\int_0^t\Vert  \naby(\eta-\zeta)\Vert_{L^2(\omega)}^2\dt'
\bigg)^{1/2}
\bigg(
\int_0^t\Vert \partial_{t'} \naby(\eta-\zeta)\Vert_{L^2(\omega)}^2\dt'
\bigg)^{1/2}.
%\\
%&\lesssim 
%\int_0^t\Vert   \eta-\zeta \Vert_{W^{2,2}(\omega)}^2\dt'
%+
%\int_0^t\Vert \partial_{t'} \naby(\eta-\zeta)\Vert_{L^2(\omega)}^2\dt'.
\end{align*}
Thus, it follows that
\begin{align*}
J_{11}
\leq
\frac{\dd}{\dt}\bigg(
&\delta
\Vert \overline{\bu} - \bv\Vert_{L^2(\Ozeta)}^2
+
\delta
\Vert  \naby^2(\eta-\zeta)\Vert_{L^2(\omega)}^{2}
+
\frac{\gamma}{2}
\int_0^t\Vert \partial_{t'} \naby(\eta-\zeta)\Vert_{L^2(\omega)}^2\dt'
\\&+c(\delta)c(\gamma)
\Vert \nabx \overline{\bu}
\Vert_{L^2(\Ozeta)}^6
\int_0^t\Vert   \eta-\zeta \Vert_{W^{2,2}(\omega)}^2\dt'
\bigg).
\end{align*}
where $\gamma>0$ is the viscoelastic coefficient in \eqref{shellEQ}.
Lastly, we use the properties of the Bogovskij operator and the  identity $\divx \overline{\bu} =(\mathbb{I}-\mathbb{B}_{\eta-\zeta}^\top):\nabx \overline{\bu}$ to obtain
\begin{align*}
\Vert  \partial_t\mathrm{Bog}_{\zeta}(\divx \overline{\bu})
\Vert_{L^{4/3}(\Ozeta)}
&\lesssim
\Vert  \partial_t(\mathbb{I}-\mathbb{B}_{\eta-\zeta}^\top):\nabx \overline{\bu}
\Vert_{W^{-1,4/3}(\Ozeta)}
+
\Vert  (\mathbb{I}-\mathbb{B}_{\eta-\zeta}^\top):\partial_t\nabx \overline{\bu}
\Vert_{W^{-1,4/3}(\Ozeta)}
\\
&\lesssim
\Vert  \partial_t(\eta-\zeta)
\Vert_{L^{2}(\omega)}
\Vert  \overline{\bu}
\Vert_{W^{2,2}(\Ozeta)}
+
\Vert  \naby^2(\eta-\zeta)
\Vert_{L^{2}(\omega)}
\Vert \partial_t \overline{\bu}
\Vert_{L^{2}(\Ozeta)}
\end{align*} 
so that
\begin{align*}
J_{12}
&\lesssim
\Vert \overline{\bu} - \bv\Vert_{L^2(\Ozeta)}^{1/2}
\Vert \nabx(\overline{\bu} - \bv)\Vert_{L^2(\Ozeta)}^{1/2}
\Vert  \partial_t(\eta-\zeta)
\Vert_{L^{2}(\omega)}
\Vert  \overline{\bu}
\Vert_{W^{2,2}(\Ozeta)}
\\
&+
\Vert \overline{\bu} - \bv\Vert_{L^2(\Ozeta)}^{1/2}
\Vert \nabx(\overline{\bu} - \bv)\Vert_{L^2(\Ozeta)}^{1/2}
\Vert  \naby^2(\eta-\zeta)
\Vert_{L^{2}(\omega)}
\Vert \partial_t \overline{\bu}
\Vert_{L^{2}(\Ozeta)}
\\
&\leq
\delta
\Vert \nabx(\overline{\bu} - \bv)\Vert_{L^2(\Ozeta)}^2
+
c(\delta)
\Vert  \partial_t(\eta-\zeta)
\Vert_{L^{2}(\omega)}^2
\Vert  \overline{\bu}
\Vert_{W^{2,2}(\Ozeta)}^2
\\
&+
c(\delta)
\Vert \overline{\bu} - \bv\Vert_{L^2(\Ozeta)}^2
+
c(\delta)
\Vert  \naby^2(\eta-\zeta)
\Vert_{L^{2}(\omega)}^2
\Vert \partial_t \overline{\bu}
\Vert_{L^{2}(\Ozeta)}^2.
\end{align*} 
If we now collect all the $J_i$'s and integrate them in time, we first observe that
\begin{equation}
\begin{aligned}
\label{10j}
\int_0^t\sum_{i=1}^{10}J_i\,\dd t'
\leq&
\big(\delta_*+\tfrac{1}{2}\big)
\int_0^t
\Vert  \overline{\bT} - \bU\Vert_{L^2(\Ozeta)}^2\dt'
+
\big(\delta+\tfrac{1}{2}\big)
\int_0^t
\Vert  \nabx(\overline{\bu} - \bv)\Vert_{L^2(\Ozeta)}^2\dt'
\\&
+
c(\delta)
\int_0^t
\Vert \overline{\bu} - \bv\Vert_{L^2(\Ozeta)}^2\dt'
+
c(\delta)
\int_0^t 
  \Vert \partial_{t'}(\eta-\zeta)\Vert_{L^{2}(\omega)}^2
   \Vert \overline{\bu} \Vert_{W^{2,2}(\Ozeta)}^2 
 \dt'
\\& 
+
c(\delta,\delta_*)\int_0^t
\Vert \eta-\zeta\Vert_{W^{2,2}(\omega)}^2
\big( 
\Vert \overline{\bu}  \Vert_{W^{2,2}(\Ozeta)}^2
+
\Vert \partial_{t'}\overline{\bu}  \Vert_{L^{2}(\Ozeta)}^2 
\big)\dt'
\\& 
+
c(\delta)\int_0^t
\Vert \eta-\zeta\Vert_{W^{2,2}(\omega)}^2
\big( 
\Vert \overline{p}  \Vert_{W^{1,2}(\Ozeta)}^2
+\Vert \overline{\bT}  \Vert_{W^{1,2}(\Ozeta)}^2
\big)\dt'.
\end{aligned}
\end{equation}
For $J_{11}$, we first note that $\Vert \nabx \overline{\bu}(t)
\Vert_{L^2(\Ozeta)}^6$ is uniformly bounded in time. This follows from the embedding 
\begin{align*}
W^{1,2}(I;L^2(\Ozeta)) \cap
L^{2}(I;W^{2,2}(\Ozeta)) \hookrightarrow
L^{\infty}(I;W^{1,2}(\Ozeta)).
\end{align*}
As such, it follows that 
\begin{equation}
\begin{aligned}
\label{11j}
\int_0^t J_{11}\,\dd t'
\leq& 
\delta \Big(
\Vert( \overline{\bu} - \bv)(t)\Vert_{L^2(\Ozeta)}^2 
+
\Vert  \naby^2(\eta-\zeta)(t) \Vert_{L^{2}(\omega )}^2
\Big)
\\&-
\delta \Big(
\Vert \overline{\bu}_0 - \bv_0\Vert_{L^2(\Ozeta)}^2 
+
\Vert  \naby^2(\eta_0-\zeta_0)  \Vert_{L^{2}(\omega )}^2
\Big)
\\&
+
\frac{\gamma}{2}
\int_0^t 
  \Vert \partial_{t'}\partial_y((\eta-\zeta)\Vert_{L^{2}(\omega)}^2
 \dt' 
+
c(\delta)c(\gamma)\int_0^t
\Vert \eta-\zeta \Vert_{W^{2,2}(\omega)}^2
\dt'.
\end{aligned}
\end{equation}
Lastly, we also obtain
\begin{equation}
\begin{aligned}
\label{12j}
\int_0^t J_{12}\,\dd t'
\leq& 
\delta
\int_0^t
\Vert\nabx( \overline{\bu} - \bv)\Vert_{L^2(\Ozeta)}^2\dt'
+
c(\delta)
\int_0^t
\Vert \overline{\bu} - \bv\Vert_{L^2(\Ozeta)}^2\dt'
\\&+
c(\delta)
\int_0^t 
  \Vert \partial_{t'}(\eta-\zeta)\Vert_{L^{2}(\omega)}^2
   \Vert \overline{\bu} \Vert_{W^{2,2}(\Ozeta)}^2 
 \dt'
\\& 
+
c(\delta)\int_0^t
\Vert \eta-\zeta\Vert_{W^{2,2}(\omega)}^2
\Vert \partial_{t'}\overline{\bu}  \Vert_{L^{2}(\Ozeta)}^2 
 \dt'.
\end{aligned}
\end{equation}
Consequently, if we integrate \eqref{velShellDiff1} in time, we can choose $\delta>0$ small enough in \eqref{10j}-\eqref{12j} to obtain 
\begin{align*}
\big(&\Vert (\overline{\bu}-\bv)(t)\Vert_{L^2(\Omega_{\zeta})}^2
+
 \Vert  \partial_t(\eta-\zeta)(t) \Vert_{L^{2}(\omega )}^2
+
\Vert  \naby^2(\eta-\zeta)(t) \Vert_{L^{2}(\omega )}^2 \big)
\\&\qquad+ 
\int_0^t
\Vert \nabx(\overline{\bu}-\bv)\Vert_{L^2(\Omega_{\zeta})}^2\dt'
+
\gamma 
\int_0^t\Vert  \partial_{t'}\naby(\eta-\zeta) \Vert_{L^{2}(\omega )}^2
\dt'
\\
\lesssim &
% \int_0^t\Vert  \partial_{t'}\naby \zeta \Vert_{L^{2}(\omega )}^2\dt'+
\Vert \overline{\bu}_0-\bv_0\Vert_{L^2(\Omega_{\zeta})}^2
+
 \Vert   \eta_\star-\zeta_\star \Vert_{L^{2}(\omega )}^2
+
\Vert  \naby^2(\eta_0-\zeta_0) \Vert_{L^{2}(\omega )}^2 
   \\&
+
\big(\delta_*+\tfrac{1}{2}\big)
\int_0^t
\Vert  \overline{\bT} - \bU\Vert_{L^2(\Ozeta)}^2\dt'
+
c 
\int_0^t
\Vert \overline{\bu} - \bv\Vert_{L^2(\Ozeta)}^2\dt'
\\&
+
c 
\int_0^t 
  \Vert \partial_{t'}(\eta-\zeta)\Vert_{L^{2}(\omega)}^2
   \Vert \overline{\bu} \Vert_{W^{2,2}(\Ozeta)}^2  
 \dt'
\\& 
+
c(\delta_*) \int_0^t
\Vert \eta-\zeta\Vert_{W^{2,2}(\omega)}^2
\big( c(\gamma)+\Vert \partial_{t'}\overline{\bu}  \Vert_{L^2(\Omega_{\zeta})}^2
+
\Vert \overline{\bu}  \Vert_{W^{2,2}(\Ozeta)}^2 
\big)\dt' 
\\& 
+
c \int_0^t
\Vert \eta-\zeta\Vert_{W^{2,2}(\omega)}^2
\big( \Vert \overline{p}  \Vert_{W^{1,2}(\Ozeta)}^2
+\Vert \overline{\bT}  \Vert_{W^{1,2}(\Ozeta)}^2
\big)\dt' 
\end{align*}
for any $t\in I$. By applying Gr\"onwall's lemma to the above (keeping \eqref{equiNorm} in mind), we obtain
%If we now substitute the estimates for the $J_i$'s back into \eqref{velShellDiff1}, integrate the resulting inequality in time, and apply Gr\"onwall's lemma, we obtain
\begin{equation}
\begin{aligned}
\label{velShellDiff2}
\big(&\Vert (\overline{\bu}-\bv)(t)\Vert_{L^2(\Omega_{\zeta})}^2
+
 \Vert  \partial_t(\eta-\zeta)(t) \Vert_{L^{2}(\omega )}^2
+
\Vert  \naby^2(\eta-\zeta)(t) \Vert_{L^{2}(\omega )}^2 \big)
\\&\qquad+ 
\int_0^t
\Vert \nabx(\overline{\bu}-\bv)\Vert_{L^2(\Omega_{\zeta})}^2\dt'
+
\gamma 
\int_0^t\Vert  \partial_{t'}\naby(\eta-\zeta) \Vert_{L^{2}(\omega )}^2
\dt'
\\
\lesssim& 
e^{\left(
c(\delta_*) \int_0^t
 \left( \Vert \partial_{t'}\overline{\bu}  \Vert_{L^2(\Omega_{\zeta})}^2
+
\Vert \overline{\bu}  \Vert_{W^{2,2}(\Ozeta)}^2
+\Vert \overline{p}  \Vert_{W^{1,2}(\Ozeta)}^2
+\Vert \overline{\bT}  \Vert_{W^{1,2}(\Ozeta)}^2
\right)\dt' \right)
}
\\&\times 
e^{\left(
c(\delta_*) \int_0^t
 \left(1+c(\gamma) 
\right)\dt' \right)
}
\bigg[
\Vert \overline{\bu}_0-\bv_0\Vert_{L^2(\Omega_{\zeta_0})}^2
+
 \Vert   \eta_\star-\zeta_\star \Vert_{L^{2}(\omega )}^2
\\&
+
\Vert  \naby^2(\eta_0-\zeta_0) \Vert_{L^{2}(\omega )}^2 
+
\big(\delta_*+\tfrac{1}{2}\big)
\int_0^t
\Vert  \overline{\bT} - \bU\Vert_{L^2(\Ozeta)}^2\dt' 
\bigg]
\end{aligned}
\end{equation}
for any $t\in I$ and any $\delta_*>0$. Recall that the constant $c(\gamma)\sim 1/\gamma$ depending on $\gamma>0$ is due to the estimate for $J_{11}$ above. 
%We now consider the equation for the difference $\overline{\bu}-\bv$ which satisfies
%\begin{align*} 
%\partial_t (\overline{\bu}-\bv)   = \Delx(\overline{\bu}-\bv) 
%-
%\nabx(\overline{p}-\pi) + \divx(\overline{\bT} -\bU)
%+ (\bv\cdot \nabx)\bv
%+
% \mathbf{h}_{\eta-\zeta}(\overline{\bu} )
%-
%\divx   \mathbb{G}_{\eta-\zeta}(\overline{\bu} ,\overline{p} ,\overline{\bT} ).
%\end{align*}
%Unfortunately, we can't test it with $\overline{\bu}-\bv$ since $\overline{\bu}$ id not solenoidal. 

\subsection{Estimate for the solute}
\label{sec:Estimate for the solute}
Having obtained the desired estimate \eqref{velShellDiff2} for the fluid-structure subproblem in the previous subsection, our goal in this subsection is to obtain an estimate for  
\begin{align*}
 \Vert (\overline{\rho}-\sigma)(t)\Vert_{L^2(\Omega_{\zeta})}^2
+
\Vert (\overline{\bT}-\bU)(t)\Vert_{L^2(\Omega_{\zeta})}^2 
\end{align*}
for any $t\in I$. 
We begin with the first term in the summand. By testing \eqref{rhoEquL} with $\sigma$, testing \eqref{rhoEqu2} with $\overline{\rho}$, and applying Reynold's transport theorem, we obtain the  identities
\begin{align*}
\frac{1}{2}\frac{\dd}{\dt} \Vert  \sigma(t) \Vert_{L^{2}( \Omega_{\zeta})}^2
=&0,
\\
\frac{1}{2}\frac{\dd}{\dt} \Vert  \overline{\rho}(t) \Vert_{L^{2}( \Omega_{\zeta})}^2 
=&
\frac{1}{2}
\int_{\partial \Omega_{\zeta}}(\bn\partial_t\zeta)\circ\bm{\varphi}_\zeta^{-1}\cdot\bn_\zeta\vert\overline{\rho}\vert^2\dd\mathcal{H}^1
-
\int_{\Omega_{\zeta}}\overline{\rho}(\overline{\bu}\cdot\nabx)\overline{\rho}\dx
\\
&+
\int_{\Omega_{\zeta}}\overline{\rho}[\varepsilon\Delx \overline{\rho} 
-
\varepsilon
\divx(\nabx \overline{\rho} (\mathbb{I}-\mathbb{A}_{\eta-\zeta}))
+
h_{\eta-\zeta}(\overline{\rho} ,\overline{\bu} )]\dx,
\\
- \frac{\dd}{\dt}\int_{\Omega_{\zeta}}\overline{\rho}\sigma\dx
=&
-
\int_{\partial \Omega_{\zeta}}(\bn\partial_t\zeta)\circ\bm{\varphi}_\zeta^{-1}\cdot\bn_\zeta\overline{\rho}\sigma\dd\mathcal{H}^1
-
\int_{\Omega_{\zeta}}\overline{\rho}\partial_t\sigma\dx
-
\int_{\Omega_{\zeta}}\sigma\partial_t\overline{\rho}\dx
\\
=&
-
\int_{\partial \Omega_{\zeta}}(\bn\partial_t\zeta)\circ\bm{\varphi}_\zeta^{-1}\cdot\bn_\zeta\overline{\rho}\sigma\dd\mathcal{H}^1
+
\int_{\Omega_{\zeta}}\overline{\rho} (\bv\cdot \nabx) \sigma\dx
\\&
+
\int_{\Omega_{\zeta}} \sigma(\overline{\bu} \cdot \nabx) \overline{\rho} \dx
-
\int_{\Omega_{\zeta}}\sigma\,\varepsilon\Delx \overline{\rho}  \dx
\\&
+
\int_{\Omega_{\zeta}}\sigma[
\varepsilon
\divx(\nabx \overline{\rho} (\mathbb{I}-\mathbb{A}_{\eta-\zeta}))
-
h_{\eta-\zeta}(\overline{\rho} ,\overline{\bu} )]\dx
.
\end{align*}
Thus, it follows
% from the basic identity $\tfrac{1}{2}\vert a -b\vert^2= \tfrac{1}{2}\vert a \vert^2+\tfrac{1}{2}\vert  b\vert^2-ab$ 
that
\begin{align*}
\frac{1}{2}\frac{\dd}{\dt} \Vert  (\overline{\rho}-\sigma)(t) \Vert_{L^{2}( \Omega_{\zeta})}^2 
=&
\int_{\Omega_{\zeta}}[\overline{\rho} (\bv\cdot \nabx) \sigma
+
 \sigma(\overline{\bu} \cdot \nabx) \overline{\rho}
 -
 \overline{\rho}(\overline{\bu} \cdot \nabx) \overline{\rho}] \dx
\\&+
\frac{1}{2}
\int_{\partial \Omega_{\zeta}}(\bn\partial_t\zeta)\circ\bm{\varphi}_\zeta^{-1}\cdot\bn_\zeta (\vert\overline{\rho}\vert^2-2\sigma\overline{\rho}) \dd\mathcal{H}^1
\\&
-
\int_{\Omega_{\zeta}}(\overline{\rho}-\sigma)
\varepsilon
\divx(\nabx \overline{\rho} (\mathbb{I}-\mathbb{A}_{\eta-\zeta}))\dx
\\
&+\int_{\Omega_{\zeta}}(\overline{\rho}-\sigma)\varepsilon\Delx \overline{\rho} \dx
+
\int_{\Omega_{\zeta}}(\overline{\rho}-\sigma)
h_{\eta-\zeta}(\overline{\rho} ,\overline{\bu} )\dx
\\&=:H_1+\ldots+H_5.
\end{align*}
Since using \eqref{divfreeL}, \eqref{interfaceL} and the divergence theorem leads to
\begin{align*}
\frac{1}{2}
\int_{\partial \Omega_{\zeta}}(\bn\partial_t\zeta)\circ\bm{\varphi}_\zeta^{-1}\cdot\bn_\zeta (\vert\overline{\rho}\vert^2-2\sigma\overline{\rho}) \dd\mathcal{H}^1
&=
%\frac{1}{2}
%\int_{ \Omega_{\zeta}} \divx(\bv (\vert\overline{\rho}\vert^2-2\sigma\overline{\rho})) \dx
%=
\frac{1}{2}
\int_{ \Omega_{\zeta}} (\bv\cdot\nabx) (\vert\overline{\rho}\vert^2-2\sigma\overline{\rho}) \dx
%\\&=
%\int_{ \Omega_{\zeta}}[\overline{\rho} (\bv\cdot\nabx) \overline{\rho}-(\bv\cdot\nabx)(\sigma\overline{\rho}) ]\dx
\\&=
\int_{ \Omega_{\zeta}}[\overline{\rho} (\bv\cdot\nabx) \overline{\rho}-\overline{\rho} (\bv\cdot\nabx)\sigma
-
\sigma(\bv\cdot\nabx)\overline{\rho}] \dx,
\end{align*}
it follows from Ladyzhenskaya's inequality and the fact that $\Vert \nabx \overline{\rho}(t) \Vert_{L^2(\Omega_{\zeta})}$ is essentially bounded in time that
\begin{equation}
\begin{aligned}
\label{j1j2}
\int_0^t
(H_1+H_2)\dt' 
&=
\int_0^t
\int_{\Omega_{\zeta}}
 (\sigma-\overline{\rho})((\overline{\bu} - \bv )\cdot \nabx) \overline{\rho} \dx\dt'
 \\
 &\lesssim
\int_0^t\big(
\Vert \overline{\rho}-\sigma\Vert_{L^2(\Omega_{\zeta})} \Vert \overline{\bu} - \bv \Vert_{L^2(\Omega_{\zeta})}^{1/2}
\Vert \nabx( \overline{\bu} - \bv) \Vert_{L^2(\Omega_{\zeta})}^{1/2}
\\&\quad
\times
\Vert \nabx \overline{\rho} \Vert_{L^2(\Omega_{\zeta})}^{1/2}
\Vert  \overline{\rho} \Vert_{W^{2,2}(\Omega_{\zeta})}^{1/2}\big)\dt'
\\
&\leq
\delta
\sup_{t'\in (0,t]}\Vert (\overline{\rho}-\sigma)(t')\Vert_{L^2(\Omega_{\zeta})}^2
+
\delta
\int_0^t 
\Vert \nabx( \overline{\bu} - \bv) \Vert_{L^2(\Omega_{\zeta})}^2\dt'
\\&\quad
+
c(\delta)
\int_0^t
 \Vert \overline{\bu} - \bv \Vert_{L^2(\Omega_{\zeta})}^2
\Vert  \overline{\rho} \Vert_{W^{2,2}(\Omega_{\zeta})}^2\dt'
\end{aligned}
\end{equation}
holds for any $\delta>0$ and for any $t\in I$. Next, since
\begin{align*}
\int_0^t
\Vert  \overline{\rho}-\sigma \Vert_{L^{\infty}( \Omega_{\zeta})}^2
\Vert \varepsilon \overline{\rho} \Vert_{W^{1,2}( \Omega_{\zeta})}^2 \dt'
\lesssim
\varepsilon^2
\sup_{t'\in (0,t]}\Vert \overline{\rho}(t') \Vert_{W^{1,2}( \Omega_{\zeta})}^2
\int_0^t\big(
\Vert \overline{\rho} \Vert_{W^{2,2}( \omega)}^2 
+
\Vert  \sigma\Vert_{L^{\infty}( \omega)}^2 
\big)
\dt'
\end{align*}
holds and the norms on the right are finite, we also have that 
\begin{equation}
\begin{aligned}
\label{j3}
\int_0^t
H_3\dt'
\lesssim&
\int_0^t
\Vert  \overline{\rho}-\sigma \Vert_{L^{2}( \Omega_{\zeta})}
\Vert \varepsilon \overline{\rho} \Vert_{W^{2,2}( \Omega_{\zeta})}
\Vert \eta-\zeta \Vert_{W^{1,\infty}( \omega)} \dt'
\\
&+
\int_0^t
\Vert  \overline{\rho}-\sigma \Vert_{L^{\infty}( \Omega_{\zeta})}
\Vert \varepsilon \overline{\rho} \Vert_{W^{1,2}( \Omega_{\zeta})}
\Vert \eta-\zeta \Vert_{W^{2,2}( \omega)} \dt'
%\\
%\leq&
%\delta
%\sup_{t'\in (0,t]}\Vert (\overline{\rho}-\sigma)(t')\Vert_{L^2(\Omega_{\zeta})}^2
%+
%c(\delta)\varepsilon^2
%\int_0^t
%\Vert \eta-\zeta \Vert_{W^{2,2}( \omega)}^2 
%\Vert \overline{\rho} \Vert_{W^{2,2}( \Omega_{\zeta})}^2
%\dt'
%\\&
%+
%\delta
%\sup_{t'\in (0,t]}
%\Vert (\eta-\zeta)(t') \Vert_{W^{2,2}( \omega)}^2 
%+
%c(\delta)\varepsilon^2
%\sup_{t'\in (0,t]}\Vert \overline{\rho}(t') \Vert_{W^{1,2}( \Omega_{\zeta})}^2
%\int_0^t\big(
%\Vert \overline{\rho} \Vert_{W^{2,2}( \omega)}^2 
%+
%\Vert  \sigma\Vert_{L^{\infty}( \omega)}^2 
%\big)
%\dt'
\\
\leq&
\delta
\sup_{t'\in (0,t]}\Vert (\overline{\rho}-\sigma)(t')\Vert_{L^2(\Omega_{\zeta})}^2
+
c(\delta)\varepsilon^2
\int_0^t
\Vert \eta-\zeta \Vert_{W^{2,2}( \omega)}^2 
\Vert \overline{\rho} \Vert_{W^{2,2}( \Omega_{\zeta})}^2
\dt'
\\&
+
\delta
\sup_{t'\in (0,t]}
\Vert (\eta-\zeta)(t') \Vert_{W^{2,2}( \omega)}^2 
+
c\,\varepsilon^2 
\end{aligned} 
\end{equation}
and
\begin{equation}
\begin{aligned}
\label{j4}
\int_0^t
H_4\dt' 
&\leq
\delta
\sup_{t'\in (0,t]}
\Vert (\overline{\rho}-\sigma)(t) \Vert_{L^{2}( \Omega_{\zeta})}^2
+
c(\delta)\varepsilon^2 
\int_0^t
\Vert \overline{\rho} \Vert_{W^{2,2}( \Omega_{\zeta})}^2
\dt'
\\
&\leq
\delta
\sup_{t'\in (0,t]}
\Vert (\overline{\rho}-\sigma)(t') \Vert_{L^{2}( \Omega_{\zeta})}^2
+
c(\delta) \varepsilon^2
\end{aligned}
\end{equation}
holds for any $\delta>0$ and for any $t\in I$. For $H_5$, we use the definition of $h_{\eta-\zeta}(\overline{\rho} ,\overline{\bu} )$ given after \eqref{solute2} to rewrite it as 
\begin{align*}
\int_0^t
H_5\dt' 
&
=
\int_0^t
\int_{\Omega_{\zeta}}(\overline{\rho}-\sigma)
 (1-J_{\eta-\zeta})\partial_{t'}\overline{\rho} 
 \dx\dt' 
-
\int_0^t
\int_{\Omega_{\zeta}}(\overline{\rho}-\sigma)
 J_{\eta-\zeta} \nabx \overline{\rho} \cdot\partial_{t'}\bm{\Psi}_{\eta-\zeta}^{-1}\circ \bm{\Psi}_{\eta-\zeta}  \dx\dt' 
\\&\qquad 
+
\int_0^t
\int_{\Omega_{\zeta}}(\overline{\rho}-\sigma) 
\overline{\bu}  \cdot\nabx \overline{\rho}  (\mathbb{I}-\mathbb{B}_{\eta-\zeta})]\dx\dt' 
\\&
=:\int_0^tH_5^a\dt'+\int_0^tH_5^b\dt'+\int_0^tH_5^c\dt'
\end{align*} 
where
\begin{align*}
\int_0^t
H_5^a\dt' 
&\lesssim
\int_0^t
\Vert
\overline{\rho}-\sigma\Vert_{L^2(\Ozeta)}
\Vert \eta-\zeta \Vert_{W^{1,\infty}(\omega)}\Vert\partial_{t'}\overline{\rho} \Vert_{L^2(\Ozeta)}\dt'
\\
&\leq
\delta
\sup_{t'\in (0,t]}
\Vert (\overline{\rho}-\sigma)(t') \Vert_{L^{2}( \Omega_{\zeta})}^2
+
c(\delta) 
\int_0^t
\Vert \eta-\zeta \Vert_{W^{2,2}(\omega)}^2
\Vert \partial_{t'}\overline{\rho} \Vert_{L^{2}( \Omega_{\zeta})}^2
\dt'.
\end{align*}
On the other hand, since
\begin{align*}
\sup_{t'\in (0,t]}\Vert
(\overline{\rho}-\sigma)(t')\Vert_{L^6(\Ozeta)}^2
\lesssim \sup_{t'\in (0,t]}
\big(
\Vert  \overline{\rho}(t')  \Vert_{W^{1,2}( \Omega_{\zeta})}^2
+
\Vert \sigma(t') \Vert_{L^{\infty}( \Omega_{\zeta})}^2
\big)
\end{align*}
is finite, it follows that
\begin{align*}
\int_0^t
H_5^b\dt' 
&\lesssim
\int_0^t
\Vert
\overline{\rho}-\sigma\Vert_{L^6(\Ozeta)}
\Vert \eta-\zeta \Vert_{W^{1,6}(\omega)}\Vert \overline{\rho} \Vert_{W^{1,6}(\Ozeta)}
\Vert \partial_{t'}(\eta-\zeta) \Vert_{L^{2}(\omega)}\dt'
%\\
%&\leq
%\delta
%\sup_{t'\in (0,t]}
%\Vert \partial_{t'}(\eta-\zeta)(t') \Vert_{L^{2}(\omega)}^2
%+
%c(\delta,t) \sup_{t'\in (0,t]}
%\big(
%\Vert  \overline{\rho}(t')  \Vert_{W^{1,2}( \Omega_{\zeta})}^2
%+
%\Vert \sigma(t') \Vert_{L^{\infty}( \Omega_{\zeta})}^2
%\big)
%\int_0^t
%\Vert \eta-\zeta \Vert_{W^{2,2}(\omega)}^2
%\Vert \overline{\rho} \Vert_{W^{2,2}( \Omega_{\zeta})}^2
%\dt'
\\
&\leq
\delta
\sup_{t'\in (0,t]}
\Vert \partial_{t'}(\eta-\zeta)(t') \Vert_{L^{2}(\omega)}^2
+
c(\delta) 
\int_0^t
\Vert \eta-\zeta \Vert_{W^{2,2}(\omega)}^2
\Vert \overline{\rho} \Vert_{W^{2,2}( \Omega_{\zeta})}^2
\dt'
\end{align*} 
and by using the fact that $\Vert  \overline{\bu}(t)  \Vert_{W^{1,2}( \Omega_{\zeta})}^2$ is essentially bounded in time, we obtain
\begin{align*}
\int_0^t
H_5^c\dt' 
&\lesssim
\int_0^t
\Vert
\overline{\rho}-\sigma\Vert_{L^2(\Ozeta)}
\Vert
\overline{\bu} \Vert_{L^4(\Ozeta)}
\Vert \overline{\rho} \Vert_{W^{1,4}(\Ozeta)}
\Vert \eta-\zeta\Vert_{W^{1,\infty}(\omega)}
\dt'
%\\
%&\leq
%\delta
%\sup_{t'\in (0,t]}
%\Vert  (\overline{\rho}-\sigma)(t') \Vert_{L^{2}(\omega)}^2
%+
%c(\delta) \sup_{t'\in (0,t]} 
%\Vert  \overline{\bu}(t')  \Vert_{W^{1,2}( \Omega_{\zeta})}^2
%\int_0^t
%\Vert \eta-\zeta \Vert_{W^{2,2}(\omega)}^2
%\Vert \overline{\rho} \Vert_{W^{2,2}( \Omega_{\zeta})}^2
%\dt'
\\
&\leq
\delta
\sup_{t'\in (0,t]}
\Vert  (\overline{\rho}-\sigma)(t') \Vert_{L^{2}(\omega)}^2
+
c(\delta) 
\int_0^t
\Vert \eta-\zeta \Vert_{W^{2,2}(\omega)}^2
\Vert \overline{\rho} \Vert_{W^{2,2}( \Omega_{\zeta})}^2
\dt'.
\end{align*}
We have thus shown that 
\begin{equation}
\begin{aligned}
\label{diffEstRho}
\Vert (\overline{\rho}-\sigma)(t)\Vert_{L^2(\Omega_{\zeta})}^2
\lesssim&
\Vert \overline{\rho}_0-\sigma_0\Vert_{L^2(\Omega_{\zeta_0})}^2
+
\delta
\int_0^t
\Vert \nabx( \overline{\bu} - \bv) \Vert_{L^2(\Omega_{\zeta})}^2\dt'
\\&+
\delta
\sup_{t'\in (0,t]}
\Vert  (\overline{\rho}-\sigma)(t') \Vert_{L^{2}(\omega)}^2
+
\delta
\sup_{t'\in (0,t]}
\Vert \partial_{t'}(\eta-\zeta)(t') \Vert_{L^{2}( \omega)}^2 
\\&+
\delta
\sup_{t'\in (0,t]}
\Vert (\eta-\zeta)(t') \Vert_{W^{2,2}( \omega)}^2 
\\&+c(\delta) 
\int_0^t
\big(
 \Vert \overline{\bu} - \bv \Vert_{L^2(\Omega_{\zeta})}^2
 +
 \Vert \eta-\zeta \Vert_{W^{2,2}(\omega)}^2
 \big)
\Vert  \overline{\rho} \Vert_{W^{2,2}(\Omega_{\zeta})}^2\dt'  
\\&
+c(\delta) 
\int_0^t
\Vert \eta-\zeta \Vert_{W^{2,2}(\omega)}^2
\Vert \partial_{t'}\overline{\rho} \Vert_{L^{2}( \Omega_{\zeta})}^2
\dt'
+
\varepsilon^2c(\delta) 
\end{aligned}
\end{equation} 
holds for any $\delta>0$ and for any $t\in I$.
We are now going to obtain an estimate for  $\Vert (\overline{\bT}-\bU)(t)\Vert_{L^2(\Omega_{\zeta})}^2$.  First of all, we note that due to Proposition \ref{prop:zeroCorotational}, if we test \eqref{soluteL} with $\bU$, we obtain
\begin{align*}
\frac{1}{2}\frac{\dd}{\dt} \Vert  \bU(t) \Vert_{L^{2}( \Omega_{\zeta})}^2 
=
-2
\int_{\Omega_{\zeta}}
(\bU- \sigma\mathbb{I}):\bU\dx.
\end{align*} 
If we also test \eqref{solute2} with $\overline{\bT}$, we obtain
\begin{align*}
\frac{1}{2}\frac{\dd}{\dt} \Vert  \overline{\bT}(t) \Vert_{L^{2}( \Omega_{\zeta})}^2 
=&
\frac{1}{2}
\int_{\partial \Omega_{\zeta}}(\bn\partial_t\zeta)\circ\bm{\varphi}_\zeta^{-1}\cdot\bn_\zeta\vert\overline{\bT}\vert^2\dd\mathcal{H}^1
-
\int_{\Omega_{\zeta}}(\overline{\bu}\cdot\nabx)\overline{\bT} :\overline{\bT}\dx
\\
&+
\int_{\Omega_{\zeta}} [\varepsilon\Delx \overline{\bT} 
-
\varepsilon
\divx(\nabx \overline{\bT} (\mathbb{I}-\mathbb{A}_{\eta-\zeta}))
+
\mathbb{H}_{\eta-\zeta}(\overline{\rho} ,\overline{\bu} )]: \overline{\bT}\dx
\\
&
 -2
\int_{\Omega_{\zeta}}
( \overline{\bT}-\overline{\rho}\mathbb{I}): \overline{\bT}\dx.
\end{align*}
Furthermore,
\begin{align*}
- \frac{\dd}{\dt}\int_{\Omega_{\zeta}}\overline{\bT} :\bU\dx
=&
-
\int_{\partial \Omega_{\zeta}}(\bn\partial_t\zeta)\circ\bm{\varphi}_\zeta^{-1}\cdot\bn_\zeta\overline{\bT} :\bU\dd\mathcal{H}^1
-
\int_{\Omega_{\zeta}} \partial_t\bU:\overline{\bT}\dx
-
\int_{\Omega_{\zeta}}\partial_t\overline{\bT} :\bU\dx
\\
=&
-
\int_{\partial \Omega_{\zeta}}(\bn\partial_t\zeta)\circ\bm{\varphi}_\zeta^{-1}\cdot\bn_\zeta\overline{\bT} :\bU \dd\mathcal{H}^1
+
\int_{\Omega_{\zeta}}( \bv\cdot \nabx) \bU:\overline{\bT} \dx
\\&
-
\int_{\Omega_{\zeta}}[\mathbb{W}(\nabx\bv) \bU
+
  \bU\mathbb{W}((\nabx\bv )^\top)
-
2
(\bU -\sigma \mathbb{I})]:\overline{\bT} \dx
\\&
+
\int_{\Omega_{\zeta}} [(\overline{\bu} \cdot \nabx) \overline{\bT}
+
2
(\overline{\bT} -\overline{\rho} \mathbb{I})]:\bU \dx
\\&-
\int_{\Omega_{\zeta}}[\mathbb{W}(\nabx\overline{\bu} )\overline{\bT} 
+
\overline{\bT} \mathbb{W}( (\nabx\overline{\bu} )^\top)
]:\bU\dx
\\&
-
\int_{\Omega_{\zeta}}[\varepsilon\Delx \overline{\bT} 
-
\varepsilon
\divx(\nabx \overline{\bT} (\mathbb{I}-\mathbb{A}_{\eta-\zeta}))
+
\mathbb{H}_{\eta-\zeta}(\overline{\rho} ,\overline{\bu} ,\overline{\bT})]:\bU\dx.
\end{align*}
We  now note that
\begin{align*}
-2(\bU- \sigma\mathbb{I}):\bU
-&
2( \overline{\bT}-\overline{\rho}\mathbb{I}): \overline{\bT}
+
2
(\bU -\sigma \mathbb{I}):\overline{\bT} 
+
2
(\overline{\bT} -\overline{\rho} \mathbb{I})]:\bU
\\&
=2(\overline{\rho}- \sigma) \mathbb{I}:(\overline{\bT}-\bU)
-
2\vert \overline{\bT}-\bU \vert^2
\\&
\leq\vert \overline{\rho}- \sigma\vert^2
-
\vert \overline{\bT}-\bU \vert^2
\end{align*}
and similar to \eqref{j1j2},
\begin{equation}
\begin{aligned}
\label{j1j2X}
\int_{\Omega_{\zeta}}&\big[(\bv\cdot \nabx) \bU:\overline{\bT} 
+
(\overline{\bu} \cdot \nabx) \overline{\bT}:\bU
-
(\overline{\bu}\cdot\nabx)\overline{\bT}:\overline{\bT}
\big]\dx
\\&+
\frac{1}{2}
\int_{\partial \Omega_{\zeta}}(\bn\partial_t\zeta)\circ\bm{\varphi}_\zeta^{-1}\cdot\bn_\zeta\big(\vert\overline{\bT}\vert^2
-2\overline{\bT}:\bU
\big)
\dd\mathcal{H}^1
\\&
=
\int_{\Ozeta}((\overline{\bu} - \bv )\cdot \nabx) \overline{\bT}:  (\bU - \overline{\bT} )\dx.
\end{aligned}
\end{equation}
If we combine all the information above, we obtain
\begin{equation}
\begin{aligned}
\label{diffSolu1}
\frac{1}{2}\frac{\dd}{\dt} \Vert  (\overline{\bT}-\bU)(t) \Vert_{L^{2}( \Omega_{\zeta})}^2 
&+
\Vert   \overline{\bT}-\bU \Vert_{L^{2}( \Omega_{\zeta})}^2 
\leq
\Vert   \overline{\rho}-\sigma \Vert_{L^{2}( \Omega_{\zeta})}^2 
\\&
+
\int_{\Ozeta}((\overline{\bu} - \bv )\cdot \nabx) \overline{\bT}:  (\bU - \overline{\bT} )\dx
\\&
 -
\int_{\Omega_{\zeta}}[\mathbb{W}(\nabx\bv) \bU
+
  \bU\mathbb{W}((\nabx\bv )^\top)
]: \overline{\bT}  \dx 
\\
&
-
\int_{\Omega_{\zeta}}[\mathbb{W}(\nabx \overline{\bu})\overline{\bT} + \overline{\bT}\mathbb{W}((\nabx \overline{\bu})^\top)]:\bU\dx
\\
&+
\int_{\Omega_{\zeta}} \varepsilon\Delx \overline{\bT} 
: (\overline{\bT}-\bU)\dx
\\&-
\int_{\Omega_{\zeta}} 
\varepsilon
\divx(\nabx \overline{\bT} (\mathbb{I}-\mathbb{A}_{\eta-\zeta}))
 : (\overline{\bT}-\bU)\dx
\\&+
\int_{\Omega_{\zeta}} 
\mathbb{H}_{\eta-\zeta}(\overline{\rho} ,\overline{\bu} ,\overline{\bT}): (\overline{\bT}-\bU)\dx
\\&=:K_1+\ldots+K_7.
\end{aligned}
\end{equation}
We now observe that by using \eqref{diffEstRho}, the estimate
\begin{equation}
\begin{aligned}
\label{diffEstRho1}
\int_0^t
K_1\dt'
&\lesssim
T\sup_{t'\in (0,t]}\Vert (\overline{\rho}-\sigma)(t')\Vert_{L^2(\Omega_{\zeta})}^2
\\ \lesssim&
T\bigg[
\Vert \overline{\rho}_0-\sigma_0\Vert_{L^2(\Omega_{\zeta_0})}^2
+
\delta
\int_0^t
\Vert \nabx( \overline{\bu} - \bv) \Vert_{L^2(\Omega_{\zeta})}^2\dt'
\\&+
\delta
\sup_{t'\in (0,t]}
\Vert \partial_{t'}(\eta-\zeta)(t') \Vert_{L^{2}( \omega)}^2 
 +
\delta
\sup_{t'\in (0,t]}
\Vert (\eta-\zeta)(t') \Vert_{W^{2,2}( \omega)}^2 
\\&+
c(\delta) 
\int_0^t
\big(
 \Vert \overline{\bu} - \bv \Vert_{L^2(\Omega_{\zeta})}^2
 +
 \Vert \eta-\zeta \Vert_{W^{2,2}(\omega)}^2
 \big)
\Vert  \overline{\rho} \Vert_{W^{2,2}(\Omega_{\zeta})}^2\dt'  
\\&
+
c(\delta) 
\int_0^t
\Vert \eta-\zeta \Vert_{W^{2,2}(\omega)}^2
\Vert \partial_{t'}\overline{\rho} \Vert_{L^{2}( \Omega_{\zeta})}^2
\dt'
+
\varepsilon^2c(\delta) 
\bigg]
\end{aligned}
\end{equation}
holds for any $\delta>0$ and for any $t\in I$. 
Similar to \eqref{j1j2}, we also have
\begin{equation}
\begin{aligned}
\label{k2}
\int_0^t
K_2\dt'  
\leq&
\delta
\sup_{t'\in (0,t]}\Vert (\overline{\bT}-\bU)(t')\Vert_{L^2(\Omega_{\zeta})}^2
+
\delta
\int_0^t
\Vert \nabx( \overline{\bu} - \bv) \Vert_{L^2(\Omega_{\zeta})}^2\dt'
\\&+
c(\delta)
\int_0^t
 \Vert \overline{\bu} - \bv \Vert_{L^2(\Omega_{\zeta})}^2
\Vert  \overline{\bT} \Vert_{W^{2,2}(\Omega_{\zeta})}^2\dt'
\end{aligned}
\end{equation}
for any $\delta>0$ and for any $t\in I$.
To estimate $K_3$ and $K_4$,  we first use the identity $\mathbb{W}((\nabx\bu )^\top)=-\mathbb{W}(\nabx\bu)$ and the relation $(\mathbb{A}\mathbb{B}):\mathbb{C}=(\mathbb{C}\mathbb{A}):\mathbb{B}$ which \footnote{Some authors define $\mathbb{A}:\mathbb{B}:=\sum a_{ij}b_{ij}$ rather than $\sum a_{ij}b_{ji}$ as stated in Section \ref{sec:prelim}. In that case, the same conclusion holds but with the identity $(\mathbb{A}\mathbb{B}):\mathbb{C}=(\mathbb{C}^\top\mathbb{A}):\mathbb{B}^\top$.}  holds for matrices $\mathbb{A},\mathbb{B},\mathbb{C}\in \mathbb{R}^{2\times2}$ to rewrite the terms in $K_3$ and $K_4$ as  
\begin{align*}
[\mathbb{W}(\nabx\bv) \bU
&+
  \bU\mathbb{W}((\nabx\bv )^\top)
]: \overline{\bT}  
+
[\mathbb{W}(\nabx \overline{\bu})\overline{\bT} + \overline{\bT}\mathbb{W}((\nabx \overline{\bu})^\top)]:\bU
%\\= &
%2 \overline{\bT} \mathbb{W}(\nabx\bv) :\bU
%+
%2
%  \bU\mathbb{W}((\nabx\bv )^\top)
%: \overline{\bT}  
%\\&+
%2\bU\mathbb{W}(\nabx \overline{\bu}):\overline{\bT} 
%+
%2 \overline{\bT}\mathbb{W}((\nabx \overline{\bu})^\top):\bU
\\= &
 \overline{\bT} \mathbb{W}(\nabx\bv) :\bU
-
  \bU\mathbb{W}(\nabx\bv )
: \overline{\bT}  
\\&+
\bU\mathbb{W}(\nabx \overline{\bu}):\overline{\bT} 
-
 \overline{\bT}\mathbb{W}(\nabx \overline{\bu}):\bU
\\= &
\bU\mathbb{W}(\nabx (\overline{\bu}-\bv)):(\overline{\bT}-\bU) 
-
 (\overline{\bT}-\bU)\mathbb{W}(\nabx (\overline{\bu}-\bv)):\bU
\\=& 
\tfrac{1}{2}
\Big\{
\bU \nabx (\overline{\bu}-\bv):(\overline{\bT}-\bU) 
-
\bU (\nabx (\overline{\bu}-\bv))^\top:(\overline{\bT}-\bU) 
 \Big\}
 \\
 \quad-&\tfrac{1}{2}\Big\{
 (\overline{\bT}-\bU) \nabx (\overline{\bu}-\bv):\bU
 -
 (\overline{\bT}-\bU) (\nabx (\overline{\bu}-\bv))^\top:\bU
 \Big\}
\end{align*} 
Thus, it follows that
\begin{align*}
\int_0^t (K_3+K_4)\dt'
\leq&
2
\int_0^t \Vert \nabx( \overline{\bu} - \bv) \Vert_{L^2(\Omega_{\zeta})} 
 \Vert \overline{\bT} -\bU \Vert_{L^2(\Omega_{\zeta})} 
\Vert  \bU \Vert_{L^{\infty}(\Omega_{\zeta})} \dt'
\\
\leq
&\frac{1}{2}
\int_0^t
\Vert \nabx( \overline{\bu} - \bv) \Vert_{L^2(\Omega_{\zeta})}^2\dt'
+ 
2
\int_0^t
 \Vert \overline{\bT} -\bU \Vert_{L^2(\Omega_{\zeta})}^2
\Vert  \bU \Vert_{L^{\infty}(\Omega_{\zeta})}^2\dt'.
\end{align*} 
Also, similar to \eqref{j4},
\begin{align*} 
\int_0^t
K_5\dt'  
&\leq
\delta
\sup_{t'\in (0,t]}
\Vert (\overline{\bT}-\bU)(t') \Vert_{L^{2}( \Omega_{\zeta})}^2
+
c(\delta) \varepsilon^2
\end{align*} 
also holds for any $\delta>0$ and for any $t\in I$.  Just as in \eqref{j3}, 
\begin{align*}  
\int_0^t
K_6\dt'
\leq&
\delta
\sup_{t'\in (0,t]}\Vert (\overline{\bT}-\bU)(t')\Vert_{L^2(\Omega_{\zeta})}^2
+
c(\delta)\varepsilon^2
\int_0^t
\Vert \eta-\zeta \Vert_{W^{2,2}( \omega)}^2 
\Vert \overline{\bT} \Vert_{W^{2,2}( \Omega_{\zeta})}^2
\dt'
\\&
+
\delta
\sup_{t'\in (0,t]}
\Vert (\eta-\zeta)(t') \Vert_{W^{2,2}( \omega)}^2 
+
c(\delta) \varepsilon^2 .
\end{align*}
If we compare how the estimates for  $H_5^a,H_5^b$ and $H_5^c$ where obtained, we observe that the estimate
\begin{align*}
\int_0^t
K_7\dt' 
&\leq
\delta
\sup_{t'\in (0,t]}
\Vert (\overline{\bT}-\bU)(t') \Vert_{L^{2}( \Omega_{\zeta})}^2
+
c(\delta) 
\int_0^t
\Vert \eta-\zeta \Vert_{W^{2,2}(\omega)}^2
\Vert \partial_t\overline{\bT} \Vert_{L^{2}( \Omega_{\zeta})}^2
\dt'
\\
&+
\delta
\sup_{t'\in (0,t]}
\Vert \partial_{t'}(\eta-\zeta)(t') \Vert_{L^{2}(\omega)}^2
+
c(\delta) 
\int_0^t
\Vert \eta-\zeta \Vert_{W^{2,2}(\omega)}^2
\Vert \overline{\bT} \Vert_{W^{2,2}( \Omega_{\zeta})}^2
\dt'
\\
&+ 2
\int_0^t\int_{\Ozeta}
\vert\overline{\bT} \vert\,\vert\mathbb{B}_{\eta-\zeta}-\mathbb{I}\vert\,\vert\nabx\overline{\bu} \vert\,\vert \overline{\bT}-\bU\vert
\dx
\dt'
\\
&+ 2
\int_0^t\int_{\Ozeta}
\vert 1-J_{\eta-\zeta}\vert\,\vert\overline{\bT} -\overline{\rho} \mathbb{I}\vert\,\vert \overline{\bT}-\bU\vert
\dx
\dt'
\end{align*}
holds for any $\delta>0$ and for any $t\in I$ where
\begin{align*}
2
\int_0^t\int_{\Ozeta}&
\vert\overline{\bT} \vert\,\vert\mathbb{B}_{\eta-\zeta}-\mathbb{I}\vert\,\vert\nabx\overline{\bu} \vert\,\vert \overline{\bT}-\bU\vert
\dx\dt'
\\&\lesssim
\int_0^t
\Vert \overline{\bT}\Vert_{L^6(\Ozeta)}
\Vert \eta-\zeta \Vert_{W^{1,6}(\omega)}
\Vert \overline{\bu}\Vert_{W^{1,6}(\Ozeta)}
\Vert \overline{\bT}-\bU\Vert_{L^2(\Ozeta)}
\dt'
\\&
\leq
\delta
\sup_{t'\in (0,t]}
\Vert (\overline{\bT}-\bU)(t') \Vert_{L^{2}( \Omega_{\zeta})}^2
+
c(\delta) 
\int_0^t
\Vert \eta-\zeta \Vert_{W^{2,2}(\omega)}^2
\Vert \overline{\bu} \Vert_{W^{2,2}( \Omega_{\zeta})}^2
\dt'
\end{align*}
since $\Vert \overline{\bT}\Vert_{L^6(\Ozeta)}$ is essentially bounded in time. Also,
\begin{align*}
2
\int_0^t\int_{\Ozeta}&
\vert 1-J_{\eta-\zeta}\vert\,\vert\overline{\bT} -\overline{\rho} \mathbb{I}\vert\,\vert \overline{\bT}-\bU\vert
\dx
\dt'
\\&\lesssim
\int_0^t
\Vert \eta-\zeta \Vert_{W^{1,\infty}(\omega)} 
\Vert \overline{\bT}-\bU\Vert_{L^2(\Ozeta)}\big(
\Vert \overline{\bT}\Vert_{L^2(\Ozeta)}+
\Vert \overline{\rho}\Vert_{L^2(\Ozeta)}\big)
\dt'
\\&
\leq
\delta
\sup_{t'\in (0,t]}
\Vert (\overline{\bT}-\bU)(t') \Vert_{L^{2}( \Omega_{\zeta})}^2
+
c(\delta) 
\int_0^t
\Vert \eta-\zeta \Vert_{W^{2,2}(\omega)}^2
\big(
\Vert \overline{\bT}\Vert_{L^2(\Ozeta)}^2+
\Vert \overline{\rho}\Vert_{L^2(\Ozeta)}^2\big)
\dt'
\end{align*}
holds for any $\delta>0$ and for any $t\in I$.
\\
If we now collect all the estimates for the $K_i$'s, we obtain that
\begin{equation}
\begin{aligned} \nonumber
\int_0^t
\sum_{i=1}^7K_i\dt'
\leq& 
2
\int_0^t
 \Vert \overline{\bT} -\bU \Vert_{L^2(\Omega_{\zeta})}^2
\Vert  \bU \Vert_{L^{\infty}(\Omega_{\zeta})}^2\dt' 
+
\delta
\sup_{t'\in (0,t]}\Vert (\overline{\bT}-\bU)(t')\Vert_{L^2(\Omega_{\zeta})}^2
\\&+
cT
\Vert \overline{\rho}_0-\sigma_0\Vert_{L^2(\Omega_{\zeta_0})}^2
+
\big[\delta(1+ T)+\tfrac{1}{2}\big]
\int_0^t
\Vert \nabx( \overline{\bu} - \bv) \Vert_{L^2(\Omega_{\zeta})}^2\dt'
\\&+
\delta (1+T)
\sup_{t'\in (0,t]}
\Vert \partial_{t'}(\eta-\zeta)(t') \Vert_{L^{2}( \omega)}^2 
 +
\delta (1+T)
\sup_{t'\in (0,t]}
\Vert (\eta-\zeta)(t') \Vert_{W^{2,2}( \omega)}^2  
\\&+
c(\delta) 
\int_0^t
\big(
 \Vert \overline{\bu} - \bv \Vert_{L^2(\Omega_{\zeta})}^2
 +
 \Vert \eta-\zeta \Vert_{W^{2,2}(\omega)}^2
 \big)
 \big(T
\Vert  \overline{\rho} \Vert_{W^{2,2}(\Omega_{\zeta})}^2+
\Vert  \overline{\bT} \Vert_{W^{2,2}(\Omega_{\zeta})}^2
\big)\dt'    
\\&
+
c(\delta) 
\int_0^t
\Vert \eta-\zeta \Vert_{W^{2,2}(\omega)}^2
\big(
T\Vert \partial_{t'}\overline{\rho} \Vert_{L^{2}( \Omega_{\zeta})}^2 
+
\Vert \partial_{t'}\overline{\bT} \Vert_{L^{2}( \Omega_{\zeta})}^2 
+
\Vert  \overline{\bu} \Vert_{W^{2,2}(\Omega_{\zeta})}^2  
\big)
\dt'  
\\
&+
c(\delta) 
\int_0^t
\Vert \eta-\zeta \Vert_{W^{2,2}(\omega)}^2
\big( (1+\varepsilon^2)
\Vert \overline{\bT} \Vert_{W^{2,2}( \Omega_{\zeta})}^2
+
\Vert \overline{\rho} \Vert_{L^{2}( \Omega_{\zeta})}^2
\big)
\dt'
+
\varepsilon^2c(\delta)(1+T)
\end{aligned}
\end{equation}
holds for any $\delta>0$ and for any $t\in I$.  
%\begin{equation}
%\begin{aligned}
%\sup_{t'\in (0,t]}\Vert (\overline{\bT}-\bU)(t')\Vert_{L^2(\Omega_{\zeta})}^2
%&+
%\int_0^t
%\Vert \overline{\bT}-\bU \Vert_{L^2(\Omega_{\zeta})}^2\dt'
%\\
%\lesssim&
%e^{\left(
%\int_0^t
%\left(  \Vert  \bU \Vert_{L^{\infty}(\Omega_{\zeta})}^2 
%+
%\Vert \nabx \overline{\bu} \Vert_{L^\infty(\Omega_{\zeta})}  \right)\dt'\right) }
%\bigg[
%\int_0^t
%\Vert \overline{\rho}-\sigma \Vert_{L^2(\Omega_{\zeta})}^2\dt'
%\\&+
%\Vert (\overline{\bT}-\bU)(0)\Vert_{L^2(\Omega_{\zeta(0)})}^2
%+
%\delta
%\int_0^t
%\Vert \nabx( \overline{\bu} - \bv) \Vert_{L^2(\Omega_{\zeta})}^2\dt'
%\\&+
% \delta
%\sup_{t'\in (0,t]}
%\Vert \partial_{t'}(\eta-\zeta)(t') \Vert_{L^{2}( \omega)}^2 
% +
% \delta
%\sup_{t'\in (0,t]}
%\Vert (\eta-\zeta)(t') \Vert_{W^{2,2}( \omega)}^2 
%\\&+
%\int_0^t
%\big(
% \Vert \overline{\bu} - \bv \Vert_{L^2(\Omega_{\zeta})}^2 
% +
%  \Vert \eta-\zeta \Vert_{W^{2,2}(\omega)}^2
% \big)
%\Vert  \overline{\bT} \Vert_{W^{2,2}(\Omega_{\zeta})}^2\dt'  
%\\&
%+
%\int_0^t
%\Vert \eta-\zeta \Vert_{W^{2,2}(\omega)}^2
%\big(
%\Vert \partial_{t'}\overline{\bT} \Vert_{L^{2}( \Omega_{\zeta})}^2 
%+
%\Vert  \overline{\bu} \Vert_{W^{2,2}(\Omega_{\zeta})}^2  
%\big)
%\dt'  
%\\&
%+
%\int_0^t
%\Vert \eta-\zeta \Vert_{W^{2,2}(\omega)}^2
%\big(  
% \Vert \overline{\bT}\Vert_{L^2(\Ozeta)}^2 
%+
% \Vert \overline{\rho}\Vert_{L^2(\Ozeta)}^2 
%\big)
%\dt'
%+
%\varepsilon^2
%\bigg].
%\end{aligned}
%\end{equation}
Thus,  we obtain from \eqref{diffSolu1} and Gr\"onwall's lemma that 
\begin{equation}
\begin{aligned}
\label{diffSolu2}
\Vert (\overline{\bT}-&\bU)(t)\Vert_{L^2(\Omega_{\zeta})}^2
+
\int_0^t
\Vert \overline{\bT}-\bU \Vert_{L^2(\Omega_{\zeta})}^2\dt'
\\
\lesssim&
e^{ c 
\int_0^t
  \Vert  \bU \Vert_{L^{\infty}(\Omega_{\zeta})}^2 
 \dt' }
\bigg[
\delta
\sup_{t'\in (0,t]}
\Vert (\overline{\bT}-\bU)(t') \Vert_{L^{2}( \Omega_{\zeta})}^2
+
\Vert \overline{\bT}_0-\bU_0\Vert_{L^2(\Omega_{\zeta_0})}^2
\\&+
T\Vert  \overline{\rho}_0-\sigma_0\Vert_{L^2(\Omega_{\zeta_0})}^2
+
\big[\delta(1+ T)+\tfrac{1}{2}\big]
\int_0^t
\Vert \nabx( \overline{\bu} - \bv) \Vert_{L^2(\Omega_{\zeta})}^2\dt'
\\&+
 \delta(1+T)
\sup_{t'\in (0,t]}
\Vert \partial_{t'}(\eta-\zeta)(t') \Vert_{L^{2}( \omega)}^2 
 +
 \delta(1+T)
\sup_{t'\in (0,t]}
\Vert (\eta-\zeta)(t') \Vert_{W^{2,2}( \omega)}^2 
\\&+
c(\delta) 
\int_0^t
\big(
 \Vert \overline{\bu} - \bv \Vert_{L^2(\Omega_{\zeta})}^2 
 +
  \Vert \eta-\zeta \Vert_{W^{2,2}(\omega)}^2
 \big)
\big(
T\Vert  \overline{\rho} \Vert_{W^{2,2}(\Omega_{\zeta})}^2
+
\Vert  \overline{\bT} \Vert_{W^{2,2}(\Omega_{\zeta})}^2\big)\dt'  
\\&
+
c(\delta) 
\int_0^t
\Vert \eta-\zeta \Vert_{W^{2,2}(\omega)}^2
\big(
T\Vert \partial_{t'}\overline{\rho} \Vert_{L^{2}( \Omega_{\zeta})}^2 
+
\Vert \partial_{t'}\overline{\bT} \Vert_{L^{2}( \Omega_{\zeta})}^2 
+
\Vert  \overline{\bu} \Vert_{W^{2,2}(\Omega_{\zeta})}^2  
\big)
\dt'  
\\&
+
c(\delta) 
\int_0^t
\Vert \eta-\zeta \Vert_{W^{2,2}(\omega)}^2
\big(  (1+\varepsilon^2)
 \Vert \overline{\bT}\Vert_{W^{2,2}(\Ozeta)}^2 
+
 \Vert \overline{\rho}\Vert_{L^2(\Ozeta)}^2 
\big)
\dt'
+
\varepsilon^2c(\delta) (1+T)
\bigg].
\end{aligned}
\end{equation} 
holds for any $\delta>0$ and for any $t\in I$. In particular, the constant in the exponential function is independent of $\delta>0$. 
\subsection{Full estimate}
In order to get our desired estimate  \eqref{contrEst0}, we sum up the estimates \eqref{velShellDiff2}, \eqref{diffEstRho} and \eqref{diffSolu2} to obtain
\begin{equation}
\begin{aligned}
&\big(\Vert (\overline{\bu}-\bv)(t)\Vert_{L^2(\Omega_{\zeta})}^2
+
 \Vert  \partial_t(\eta-\zeta)(t) \Vert_{L^{2}(\omega )}^2
+
\Vert  \naby^2(\eta-\zeta)(t) \Vert_{L^{2}(\omega )}^2 \big)
\\&\qquad+ 
\int_0^t
\Vert \nabx(\overline{\bu}-\bv)\Vert_{L^2(\Omega_{\zeta})}^2\dt'
+
\gamma 
\int_0^t\Vert  \partial_{t'}\naby(\eta-\zeta) \Vert_{L^{2}(\omega )}^2
\dt' 
\\&
\qquad+\big(\Vert (\overline{\rho}-\sigma)(t'\Vert_{L^2(\Omega_{\zeta})}^2
+
\Vert (\overline{\bT}-\bU)(t)\Vert_{L^2(\Omega_{\zeta})}^2
\big)
+ 
\int_0^t
\Vert \overline{\bT}-\bU \Vert_{L^2(\Omega_{\zeta})}^2\dt'
\\
\lesssim&
e^{ 
c(\delta_*) \int_0^t
 \left( 1+c(\gamma)
+
\Vert \partial_{t'}\overline{\bu}  \Vert_{L^2(\Omega_{\zeta})}^2
+
\Vert \overline{\bu}  \Vert_{W^{2,2}(\Ozeta)}^2
+\Vert \overline{p}  \Vert_{W^{1,2}(\Ozeta)}^2
+\Vert \overline{\bT}  \Vert_{W^{1,2}(\Ozeta)}^2
+
\Vert  \bU \Vert_{L^{\infty}(\Omega_{\zeta})}^2 
\right)\dt' }
\\&\times
\bigg[ 
\Vert \overline{\bu}_0-\bv_0\Vert_{L^2(\Omega_{\zeta_0})}^2
+
 \Vert   \eta_\star-\zeta_\star \Vert_{L^{2}(\omega )}^2
+
\Vert  \naby^2(\eta_0-\zeta_0) \Vert_{L^{2}(\omega )}^2  
+
\Vert  \overline{\bT}_0 -\bU_0\Vert_{L^2(\Omega_{\zeta_0})}^2
\\&+
(1+T)\Vert  \overline{\rho}_0 - \sigma_0\Vert_{L^2(\Omega_{\zeta_0})}^2
+
\big[\delta(1+ T)+\tfrac{1}{2}\big]
\int_0^t
\Vert \nabx( \overline{\bu} - \bv) \Vert_{L^2(\Omega_{\zeta})}^2\dt'
\\&+
 \delta(1+T)
\sup_{t'\in (0,t]}
\Vert \partial_{t'}(\eta-\zeta)(t') \Vert_{L^{2}( \omega)}^2 
 +
 \delta(1+T)
\sup_{t'\in (0,t]}
\Vert (\eta-\zeta)(t') \Vert_{W^{2,2}( \omega)}^2 
\\&
+
\delta
\sup_{t'\in (0,t]}
\Vert (\overline{\rho}-\sigma)(t') \Vert_{L^{2}( \Omega_{\zeta})}^2
+
\delta
\sup_{t'\in (0,t]}
\Vert (\overline{\bT}-\bU)(t') \Vert_{L^{2}( \Omega_{\zeta})}^2 
+
\big(\delta_*+\tfrac{1}{2}\big)
\int_0^t
\Vert  \overline{\bT} - \bU\Vert_{L^2(\Ozeta)}^2\dt' 
\\&+
c(\delta) 
\int_0^t
\big(
 \Vert \overline{\bu} - \bv \Vert_{L^2(\Omega_{\zeta})}^2 
 +
  \Vert \eta-\zeta \Vert_{W^{2,2}(\omega)}^2
 \big)
\big(
(1+T)\Vert  \overline{\rho} \Vert_{W^{2,2}(\Omega_{\zeta})}^2
+
\Vert  \overline{\bT} \Vert_{W^{2,2}(\Omega_{\zeta})}^2\big)\dt'  
\\&
+
c(\delta) 
\int_0^t
\Vert \eta-\zeta \Vert_{W^{2,2}(\omega)}^2
\big(
(1+T)\Vert \partial_{t'}\overline{\rho} \Vert_{L^{2}( \Omega_{\zeta})}^2 
+
\Vert \partial_{t'}\overline{\bT} \Vert_{L^{2}( \Omega_{\zeta})}^2 
+
\Vert  \overline{\bu} \Vert_{W^{2,2}(\Omega_{\zeta})}^2  
\big)
\dt'  
\\&
+
c(\delta) 
\int_0^t
\Vert \eta-\zeta \Vert_{W^{2,2}(\omega)}^2
\big(  (1+\varepsilon^2)
 \Vert \overline{\bT}\Vert_{W^{2,2}(\Ozeta)}^2 
+
 \Vert \overline{\rho}\Vert_{L^2(\Ozeta)}^2 
\big)
\dt'
+
\varepsilon^2c(\delta) (1+T)
\bigg]
\end{aligned}
\end{equation} 
for any $\delta,\delta_*>0$ and for any $t\in I$.
We can now pick any finite $\delta_*>0$ (say $\delta=1/2$) so that
if we  choose a $\delta>0$ such that
\begin{align*}
\delta
\lesssim&
\frac{1}{1+T} 
e^{ 
c\int_0^t
 \left( 1+c(\gamma)
+
\Vert \partial_{t'}\overline{\bu}  \Vert_{L^2(\Omega_{\zeta})}^2
+
\Vert \overline{\bu}  \Vert_{W^{2,2}(\Ozeta)}^2
+\Vert \overline{p}  \Vert_{W^{1,2}(\Ozeta)}^2
+\Vert \overline{\bT}  \Vert_{W^{1,2}(\Ozeta)}^2
+
\Vert  \bU \Vert_{L^{\infty}(\Omega_{\zeta})}^2 
\right)\dt' } ,
\end{align*}
then we obtain \eqref{contrEst0} from Gr\"onwall's lemma.
\end{proof}
Having proved Proposition \ref{prop:zeroCorotational}, we now  set $(\eta, \bu, p, \rho,\bT)=(\eta^\varepsilon, \bu^\varepsilon, p^\varepsilon, \rho^\varepsilon,\bT^\varepsilon)_{\varepsilon>0}$  in  \eqref{divfree}-\eqref{solute}.  
We note that if the corresponding dataset $(\rho_0, \bT_0, \bu_0, \eta_0, \eta_\star)$ satisfies \eqref{dataConv}, then the transformed dataset $(\overline{\rho}^\varepsilon_0, \overline{\bT}^\varepsilon_0, \overline{\bu}^\varepsilon_0, \eta^\varepsilon_0, \eta^\varepsilon_\star)$ satisfies
\begin{align*}
&\eta^\varepsilon_0 \rightarrow \zeta_0 \quad\text{in}\quad W^{2,2}( \omega), 
\\&\eta^\varepsilon_\star \rightarrow \zeta_\star \quad\text{in}\quad L^2( \omega),
\\&\overline{\bu}^\varepsilon_0  \rightarrow \bv_0 \quad\text{in}\quad L^2( \Omega_{\zeta(0)}), 
\\&\overline{\rho}^\varepsilon_0 \rightarrow \sigma_0 \quad\text{in}\quad L^2( \Omega_{\zeta(0)})
, 
\\&
\overline{\bT}^\varepsilon_0  \rightarrow \bU_0 \quad\text{in}\quad L^2( \Omega_{\zeta(0)})
\end{align*}
as $\varepsilon\rightarrow0$.
Thus,  an immediate corollary of Proposition \ref{prop:zeroCorotational} is the following result.
\begin{corollary}
\label{cor:main}
Let $(\eta^\varepsilon,\overline{\bu}^\varepsilon, \overline{p}^\varepsilon,\overline{\rho}^\varepsilon,\overline{\bT}^\varepsilon)_{\varepsilon>0}$ be a collection of strong solutions of \eqref{divfree}-\eqref{interface} with dataset $(\overline{\rho}^\varepsilon_0, \overline{\bT}^\varepsilon_0, \overline{\bu}^\varepsilon_0, \eta^\varepsilon_0, \eta^\varepsilon_\star)$ for which \eqref{dataConv} holds with $(\sigma_0,\bU_0)$ satisfying \eqref{mainDataForAllBounded}. Then we have
\begin{align*}
&\eta^\varepsilon \rightarrow \zeta \quad\text{in}\quad L^\infty(I;W^{2,2}( \omega)),  
\\&\partial_t\eta^\varepsilon \rightarrow \partial_t\zeta \quad\text{in}\quad L^\infty(I;L^{2}( \omega)),
\\&\partial_t\eta^\varepsilon \rightarrow \partial_t\zeta \quad\text{in}\quad L^2(I;W^{1,2}( \omega)),
\\&\overline{\bu}^\varepsilon  \rightarrow \bv \quad\text{in}\quad L^\infty(I;L^2( \Omega_{\zeta})), 
\\
&\nabx\overline{\bu}^\varepsilon  \rightarrow \nabx\bv \quad\text{in}\quad L^2(I;L^2( \Omega_{\zeta})),
\\&\overline{\rho}^\varepsilon \rightarrow \sigma \quad\text{in}\quad L^\infty(I;L^2( \Omega_{\zeta}))
, 
\\&
\overline{\bT}^\varepsilon  \rightarrow \bU \quad\text{in}\quad L^\infty(I;L^2( \Omega_{\zeta}))
\end{align*}
where $(\zeta,\bv,  \sigma,\bU)$ is an essentially bounded weak solution of \eqref{divfreeL}-\eqref{interfaceL} with dataset $(\sigma_0, \bU_0, \bv_0, \zeta_0, \zeta_\star)$.
\end{corollary}
Now we note that since
\begin{equation}
\begin{aligned}
\label{triIneq}
\Vert \bm{1}_{\Omega_{\eta^\varepsilon}}
\bu^\varepsilon  - \bm{1}_{\Omega_{\zeta}}\bv \Vert_{L^2(\Omega \cup S_\ell)}^2
\lesssim&
\Vert(\bm{1}_{\Omega_{\eta^\varepsilon}}
  - \bm{1}_{\Omega_{\zeta}}) \bu^\varepsilon \Vert_{L^2(\Omega \cup S_\ell)}^2
+
\Vert
\bm{1}_{\Omega_{\zeta}}( \overline{\bu}^\varepsilon  - \bv)
\Vert_{L^2(\Omega \cup S_\ell)}^2
\\&+
  \Vert
\bm{1}_{\Omega_{\zeta}}(\bu^\varepsilon -\overline{\bu}^\varepsilon)
\Vert_{L^2(\Omega \cup S_\ell)}^2,
\end{aligned}
\end{equation}
the first two terms to the right converge to zero as $\varepsilon\rightarrow0$ due to Corollary \ref{cor:main}.
%
%***********************
%
%Now we note that
%\begin{align}
%\label{triIneq1}
%\Vert \bm{1}_{\Omega_{\eta^\varepsilon}}
%\bT^\varepsilon  - \bm{1}_{\Omega_{\zeta}}\bU \Vert_{L^2(\Omega \cup S_\ell)}^2
%&\lesssim
%\Vert \bm{1}_{\Omega_{\eta^\varepsilon}}  \bT^\varepsilon 
%-
%\bm{1}_{\Omega_{\zeta}} \overline{\bT}^\varepsilon
%\Vert_{L^2(\Omega \cup S_\ell)}^2
%+
%\Vert
%\bm{1}_{\Omega_{\zeta}}( \overline{\bT}^\varepsilon  - \bU)
%\Vert_{L^2(\Omega \cup S_\ell)}^2.
%\end{align}
Also, since $ \overline{\bu}^\varepsilon = \bu^\varepsilon \circ \bm{\Psi}_{\eta^\varepsilon-\zeta}= \bu^\varepsilon (t, \bx+\bn(\eta^\varepsilon-\zeta)\phi)$, the last term to the right in \eqref{triIneq} also converges to zero as $\varepsilon\rightarrow0$ due to the strong uniform convergence of $\eta^\varepsilon$ to $\zeta$. Thus, it follows that
\begin{align*}
\bm{1}_{\Omega_{\eta^\varepsilon}} \bu^\varepsilon  \rightarrow \bm{1}_{\Omega_{\zeta}}\bv \quad\text{in}\quad L^\infty(I;L^2(\Omega \cup S_\ell)).
\end{align*}
Similarly, we obtain strong convergence of $\bm{1}_{\Omega_{\eta^\varepsilon}}(\nabx\bu^\varepsilon, \rho^\varepsilon ,
\bT^\varepsilon )$ to their respective limits as stated in Theorem \ref{thm:main2} thereby finishing the proof.

\section*{Statements and Declarations} 
%\subsection*{Acknowledgements}
%The author would like to thank Dominic Breit for many useful discussions.
%\subsection*{Funding}
%This work has been supported by the European Research Council (ERC) Synergy grant STUOD-DLV-856408.
\subsection*{Author Contribution}
The author wrote and reviewed the manuscript.
\subsection*{Conflict of Interest}
The author declares that they have no conflict of interest.
\subsection*{Data Availability Statement}
Data sharing is not applicable to this article as no datasets were generated
or analyzed during the current study.
\subsection*{Competing Interests}
The author have no competing interests to declare that are relevant to the content of this article.
%****List of advanced bibliographystyle****
% 1. spbasic
% 2. spphys
% 3. spmpsci
%%
%
%\bibliographystyle{spmpsci}
%\bibliography{myBibliography}
%
%

\end{document}